\documentclass{article}
\usepackage[utf8x]{inputenc}  
\usepackage[left=3.1cm, right=3.1cm, bottom=4cm]{geometry}                 
\usepackage[textsize=small]{todonotes}

\usepackage{parskip}

\usepackage{hyperref} 
\usepackage{amssymb,amsmath}
\usepackage{bm} 
\usepackage{esint} 
\usepackage{pbox}

\usepackage{enumitem}
\setlist[enumerate,1]{label=(\roman*)}

\usepackage[backend=bibtex,style=alphabetic]{biblatex}
\numberwithin{equation}{section}
\allowdisplaybreaks

\usepackage{amsthm}

\theoremstyle{plain}
\newtheorem{theorem}{Theorem}[section]
\newtheorem{proposition}[theorem]{Proposition}
\newtheorem{corollary}[theorem]{Corollary}
\newtheorem{lemma}[theorem]{Lemma}

\newtheorem{conjecture}[theorem]{Conjecture}

\theoremstyle{definition}
\newtheorem{definition}[theorem]{Definition}
\newtheorem*{definition*}{Definition}

\theoremstyle{remark}
\newtheorem{remark}[theorem]{Remark}
\newtheorem*{remark*}{Remark}

\usepackage{zref-clever}

\AddToHook{env/definition/begin}{%
\zcsetup{countertype={theorem=definition}}}
\AddToHook{env/example/begin}{%
\zcsetup{countertype={theorem=example}}}
\AddToHook{env/lemma/begin}{%
\zcsetup{countertype={theorem=lemma}}}
\AddToHook{env/corollary/begin}{%
\zcsetup{countertype={theorem=corollary}}}
\AddToHook{env/proposition/begin}{%
\zcsetup{countertype={theorem=proposition}}}
\AddToHook{env/remark/begin}{%
\zcsetup{countertype={theorem=remark}}}
\AddToHook{env/question/begin}{%
\zcsetup{countertype={theorem=question}}}
\AddToHook{env/conclusion/begin}{%
\zcsetup{countertype={theorem=conclusion}}}
\AddToHook{env/conjecture/begin}{%
\zcsetup{countertype={theorem=conjecture}}}

\zcRefTypeSetup{item}{name-sg= ,name-pl= }
\zcRefTypeSetup{equation}{name-sg= ,name-pl= }
\zcRefTypeSetup{conjecture}{name-sg=conjecture,name-pl=conjectures,Name-sg=Conjecture,Name-pl=Conjectures}

\zcRefTypeSetup{section}{cap}
\zcRefTypeSetup{appendix}{cap}
\zcRefTypeSetup{figure}{cap}
\zcRefTypeSetup{theorem}{cap}
\zcRefTypeSetup{definition}{cap}
\zcRefTypeSetup{example}{cap}
\zcRefTypeSetup{lemma}{cap}
\zcRefTypeSetup{corollary}{cap}
\zcRefTypeSetup{proposition}{cap}
\zcRefTypeSetup{remark}{cap}
\zcRefTypeSetup{question}{cap}
\zcRefTypeSetup{conclusion}{cap}
\zcRefTypeSetup{conjecture}{cap}

\usepackage{subcaption}

\newcommand\R{\mathbb{R}}
\newcommand\N{\mathbb{N}}

\newcommand\SO[1]{\mathrm{SO}(#1)}

\newcommand{\intd}{\mathop{}\!\mathrm{d}}
\renewcommand{\d}{\intd}

\newcommand\BV{\mathrm{BV}}
\DeclareMathOperator\Div{div}
\newcommand\loc{\mathrm{loc}}
\newcommand\cpt{\mathrm{c}}

\newcommand\vt[1]{\bm{#1}}
\newcommand\x{{\vt x}}
\newcommand\y{{\vt y}}
\newcommand\z{{\vt z}}
\newcommand\vtf\vt
\newcommand\rad\widetilde
\newcommand\lowerrep\underline
\newcommand\rep[1]{{#1}^*}
\DeclareMathOperator{\var}{Var}   
\newcommand\varalt{\mathop{\widetilde{\mathrm{Var}}}}
   
\newcommand\M{\mathrm M}
\newcommand\Ms{\M^{\mathrm s}}

\newcommand\Mo{\M_{\mathrm o}}
\newcommand\Mu{\M_{\mathrm u}}
\newcommand\lmo[1]{\mathcal H^{#1}}

\newcommand\sm[1]{\lmo{d-1}(#1)}
\newcommand\grad{\mathrm D}

\newcommand\fulltitle{On the variation and asymptotics of maximal functions}

\title{\fulltitle}

\newcommand\addressictp{ICTP -- The Abdus Salam International Centre for Theoretical Physics, Strada Costiera, 11, I-34151, Trieste, Italy}
\newcommand\addressictpshort{ICTP}

\author{
Emanuel Carneiro\footnote{
\addressictp,
\texttt{carneiro@ictp.it}
},\ \ 
Julian Weigt\footnote{
\addressictpshort,
\texttt{jweigt@ictp.it}
}\ \ 
}

\date{\today} 

\begin{document}

\maketitle

\begingroup
\begin{NoHyper}%
\renewcommand\thefootnote{}\footnotetext{%
}%
\addtocounter{footnote}{-1}%
\end{NoHyper}%
\endgroup

\begin{abstract}
\noindent
In this paper we present partial progress towards two conjectures on the variation of maximal operators of convolution type $\M_{\varphi}$, associated to non-negative and radial decreasing kernels $\varphi$ of integral $1$. In dimension $d=1$ we prove that $\var(\M_{\varphi} f) \leq \var(f)$ when $f$ is a characteristic function. In higher dimensions we prove that $\var(\M_{\varphi} f) \leq C_d \var(f)$ when $f$ is a radial characteristic function, with a constant $C_d$ depending only on the dimension $d$ but not on the kernel $\varphi$. We further show that $C_d \sim \sqrt{d}\log d $ is the sharp dependence on $d$ for such an inequality. This last claim relies on a novel explicit asymptotic estimate for the Hardy-Littlewood maximal function of the unit ball in $\R^d$, a result of independent interest.
\end{abstract}

\section{Introduction}
Let $\M$ be the classical centered Hardy-Littlewood maximal operator, i.e., for $f \in L^1_{{\rm loc}}(\R^d)$, 
\begin{equation}\label{20260213_15:56}
\M f(\x) = \sup_{r>0}\dfrac{1}{|B_r|}\int_{B(\x, r)} |f(\y)|\d \y\,,
\end{equation}
where $B(\x, r)$ is the open ball centered at $\x$ with radius $r$, $B_r := B({\vt 0}, r)$, and $|B_r| = \omega_d \, r^d$ denotes its $d$-dimensional Lebesgue measure.
A celebrated result of Hardy, Littlewood and Wiener establishes that $\M: L^p(\R^d) \to L^p(\R^d)$ is bounded for $p>1$, a cornerstone in harmonic analysis. The first to consider the regularity of this operator, when acting on Sobolev functions, was Kinnunen \cite{Ki} in 1997. He showed that $\M: W^{1,p}(\R^d) \to W^{1,p}(\R^d)$ is bounded for $p>1$, and similarly for the uncentered version of this operator. The quest for analogous statements for $\BV$ and $W^{1,1}$ spaces turned out to be way more subtle. In particular, the question on whether the map $f \mapsto \grad \M f$ is bounded from $W^{1,1}(\R^d)$ to $L^1(\R^d)$, first proposed by Haj\l asz and Onninen in \cite{HO}, was settled affirmatively in dimension $d=1$ in \cite{AP, Ku, Ta} (for both the centered and uncentered cases) but is generally open in higher dimensions. A recent result of the second author \cite{Weigt1} settles this problem affirmatively in the uncentered higher dimensional case, when balls are replaced by cubes. The regularity theory of maximal operators flourished around the aforementioned question \cite{BiWe, CGR, HKS, Lu2, Weigt000, Weigt0} and some of its variants, e.g., discrete analogues \cite{BCHP, CH, LW, M}, fractional analogues \cite{BRS, CMa, KiSa, LM, Weigt00}, continuity questions \cite{BM, BMGRW, CMP, GR0, GR1, GRK0, Lu1, M2}, regularizing effects \cite{AP, GRK, LXK, WeigtLahti}, and other related themes \cite{BES, HM, KL, PPSS, Ra}. The survey paper \cite{Car} by the first author brings a more detailed account on the developments of this theory over the past decades.

In this paper we consider the more general setup of maximal functions of convolution type, motivated by the study initiated in \cite{CS, CFS}. Let $\varphi \in L^1(\R^d)$ be a non-negative function such that 
\begin{equation*}
\int_{\R^d} \varphi(\x) \d \x = 1.
\end{equation*}
For $t >0$, let $\varphi_t(\x) := t^{-d} \varphi(t^{-1}\x)$, and consider the maximal operator
\begin{equation*}
\M_{\varphi}f(\x) := \sup_{t>0}\big(|f|*\varphi_t\big)(\x).
\end{equation*}
Note that when $\varphi (\x)= \chi_{B_1}(\x)/|B_1|$, $\M_{\varphi}$ specializes to be the classical centered Hardy-Littlewood maximal operator $\M$ defined in \zcref{20260213_15:56}. 

We split our discussion here into two different parts, one corresponding to the one-dimensional setting and the other one corresponding to the higher dimensional setting. 

\subsection{The one-dimensional setting}
\label{SubSec1.1}

For an interval $I\subset\mathbb{R}$ and $f:I\rightarrow\mathbb{R}$ recall the variation of $f$,
\begin{equation}
\label{definition1dvariation}
\var_I(f)
=
\sup\left\{ \sum_{i=1}^{n} |f(x_i) - f(x_{i-1})| : n \in \mathbb{N},\ x_i\in I,\ -\infty < x_0 < x_1 < \dots < x_n < \infty \right\}
.
\end{equation}
For $I=\mathbb{R}$ we write $\var_\mathbb{R}(f)=\var(f)$.

We start by posing the following conjecture.
\begin{conjecture}\label{Conj1}
Let $\varphi:\R \to\R$ be non-negative, radial decreasing, and with $\int_{\R}\varphi(x) \d x = 1$. Then, for any $f:\R \to \R$ of bounded variation, we have
\begin{equation}\label{20260213_17:30}
\var(\M_{\varphi}f) \leq \var(f).
\end{equation}
\end{conjecture}
In other words, this \zcref*[nocap,noref]{Conj1} says that all such maximal operators are variation-diminishing. In the case of the uncentered Hardy-Littlewood maximal operator, let us call it $\widetilde{\M}$,  Aldaz and P\'{e}rez-L\'{a}zaro \cite{AP} proved that $\var(\widetilde{\M}f) \leq \var(f)$ for all $f$. The case of centered maximal operators is more difficult. For instance, for the centered Hardy-Littlewood maximal operator $\M$, in the remarkable paper \cite{Ku}, Kurka showed that 
$$\var(\M f) \leq 240004 \var(f)$$
for all $f$, but it is still unknown if such an inequality holds with constant $C=1$ (which would be sharp). As a matter of fact, in our setup, there are only two kernels $\varphi$ for which \zcref{Conj1} has been fully established \cite{CS}: the Gauss kernel $\varphi(x) = e^{-\pi |x|^2}$, and the Poisson kernel $\varphi(x) = \frac{1}{\pi (x^2 +1)} $. The proofs in \cite{CS} for these two particular kernels explore the connections with the underlying PDEs (the heat equation for the Gauss kernel, and Laplace's equation for the Poisson kernel).

In \cite{BiWe}, Bilz and the second author showed that \zcref{Conj1} is true for the centered Hardy-Littlewood maximal operator when restricted to the subclass of characteristic functions $f$; see also the related work \cite{HKS}. Our first result here extends the result of \cite{BiWe} uniformly to all kernels.

\begin{theorem}\label{thm2}
\zcref[S]{Conj1} is true for characteristic functions $f$.
\end{theorem}

In the setup of \zcref{Conj1}, a non-trivial intermediate question would be to establish the validity of \zcref{20260213_17:30} with the right-hand side multiplied by a constant $C_{\varphi}$ that may depend on $\varphi$. \zcref[S]{thm2} and the results for the Gauss and Poisson kernels in \cite{CS} provide evidence that \zcref{Conj1} should hold with the uniform constant $C=1$, which would be sharp. 

Recall that a competing weak definition of $\BV(\mathbb{R})$ and $\var(f)$ exists, which, unlike \zcref{definition1dvariation}, readily extends to $\R^d$; see \zcref{definition_bv}. While the two definitions do not coincide for generic functions when $d=1$, all inequalities in this \zcref*[nocap,noref]{SubSec1.1} hold for one definition if and only if they hold for the other, see for example \zcref{appendix_bvoned,Prop11_contbv}.

\subsection{The higher dimensional setting}
\label{subsec_higherdimensions}

In this \zcref*[noref,nocap]{subsec_higherdimensions} $\BV(\mathbb{R}^d)$ and $\var(f)$ denote the weak notions from \zcref{definition_bv}.

If $\varphi:\R^d \to\R$ is a non-negative and radial decreasing function, with $\int_{\R^d}\varphi(\x) \d \x = 1$, it is well-known (e.g., via the layer cake representation) that 
\begin{equation}\label{20260903_10:40}
\M_{\varphi}f(\x) \leq  \M f(\x)
\end{equation}
holds for all $\x \in \R^d$. In this sense, one can say that the Hardy-Littlewood maximal operator controls all the other maximal operators of convolution type (associated to radial decreasing kernels). Our next conjecture predicts a similar uniform control at the derivative level.
\begin{conjecture}\label{Conj3}
Let $\varphi:\R^d \to\R$ be non-negative, radial decreasing, and with $\int_{\R^d}\varphi(\x) \d \x = 1$. Then, for any $f \in \BV(\R^d)$, we have
\begin{equation}\label{20260213_17:41}
\var(\M_{\varphi}f) \leq C_d \var(f)\,,
\end{equation}
where $C_d$ is a constant depending only on the dimension $d$.
\end{conjecture}

One can recognize \zcref{Conj3} as a more general version of the question posed by Haj\l asz and Onninen \cite[Question 1]{HO} for the Hardy-Littlewood maximal operator. 
 
In this higher dimensional seting we are still interested in understanding in which ways an estimate like \zcref{20260213_17:41} could be sharp. A first theme of interest here is the study of the sharp dependence of $C_d$ on the dimension $d$. This is naturally motivated by several important papers in the literature that deal with dimension-free bounds (or lack thereof) for maximal operators; e.g., \cite{Aldaz, AP0, MeSo, StSt, SpectorStockdale}. As a starting point in this study, we consider the Hardy-Littlewood maximal  function of the characteristic function of the unit ball $B_1 \subset \R^d$. It is not hard to see that, for $|\x| >1$,
\begin{equation}\label{20260321_15:27}
\M\chi_{B_1}(\x) \sim_d |\x|^{-d}.
\end{equation}
In fact, this expression is ubiquitous in the literature, but we were unable to find an explicit asymptotic estimate with the correct dependence on the dimension $d$. Our next result is of independent interest and settles this problem.

\begin{theorem}\label{Thm4}
For $B_1 \subset \R^d$ and $\x \in \R^d$ with $|\x|> 1$ we have
\begin{equation}\label{20260320_15:03}
\M\chi_{B_1}(\x) \sim |\x|^{-d} \left( 1 + \frac{\sqrt{d}}{\left( |\x| + \frac{1}{\sqrt{|\x| -1}}\right)}\right)^{-1}\,,
\end{equation}
with the implicit constant being universal.
\end{theorem}
Note that an equivalent way to state \zcref{20260320_15:03} is by saying that, for $|\x| >1$, 
\begin{align}\label{20260912_16:07}
\M\chi_{B_1}(\x)
&\sim
|\x|^{-d} \times 
\begin{cases}
1\, & {\rm if} \  \ \ 1 < |\x| \leq 1 + \frac{1}{d}\,;\\ \\
\dfrac{1}{\sqrt{d}\sqrt{|\x| -1}} \, , & {\rm if} \  \ 1 + \frac{1}{d} \leq |\x| \leq 2\,;\\
\\
\dfrac{|\x|}{2\sqrt{d}} \, ,  & {\rm if} \ \ 2 \leq |\x| \leq 2\sqrt{d}\,;\\
\\
1 \, , & {\rm if} \  \ 2\sqrt{d} \leq |\x|.
\end{cases}
\end{align}


\zcref[S]{Thm4} yields several interesting corollaries. The first one is a consequence of the coarea formula, \zcref{coarea}, or more precisely \zcref{radiallydecreasingcomparison}.

\begin{corollary}\label{Cor5}
For $B_1 \subset \R^d$ we have
\begin{equation*}
\frac{\var(\M\chi_{B_1})}{\var(\chi_{B_1})} \sim \sqrt{d} \, \log d.
\end{equation*}
\end{corollary}
We remark that classical bounds in the spirit of \zcref{20260321_15:27}, like
\begin{equation}\label{20260912_17:58}
\min\left \{1, \frac{1}{(|\x|-1)^d}\right\} \geq \M\chi_{B_1}(\x) \geq \frac{1}{(|\x| +1)^d}
\end{equation}
for $|\x| >1$, only yield, via the coarea formula, 
\begin{align}\label{20260912_20:23}
2^d \gtrsim \frac{\var(\M\chi_{B_1})}{\var(\chi_{B_1})} \gtrsim 1.
\end{align}
This computation is carried out in \zcref{Remark_coarea_HL} at the end of \zcref{Sec5_Cor1.5}. In light of \zcref{20260912_17:58} and \zcref{20260912_20:23}, one can start to appreciate \zcref{Cor5} and the impact of the precise asymptotics given in \zcref[S]{Thm4}.

As a consequence of \zcref{Cor5}, we plainly conclude the following.
\begin{corollary}\label{Cor1.6_growth}
The constant $C_d$ in \zcref{Conj3} must verify
\begin{equation*}
C_d \gtrsim \sqrt{d} \, \log d.
\end{equation*}
\end{corollary}
This emphasizes a significant difference between the case $p=1$ and the case $p>1$. If $f \in W^{1,p}(\R^d)$, for $p>1$, the original argument of Kinnunen \cite{Ki} yields that $\M_{\varphi}f$ is weakly differentiable and $|\grad \M_{\varphi}f | \leq \M_{\varphi} |\grad f |$. Hence, 
\begin{align*}
\|\grad \M_{\varphi}f \|_{L^p(\mathbb{R}^d)} \leq \|\M_{\varphi} |\grad f |\|_{L^p(\mathbb{R}^d)} \leq C_p \||\grad f |\|_{L^p(\mathbb{R}^d)} \,,
\end{align*} 
with a constant $C_p$ that is dimension free due to a result of Stein \cite{Stein}. Note also the recent work of Spector and Stockdale \cite{SpectorStockdale} in which the same growth rate $\sqrt{d} \, \log d$ as in \zcref{Cor1.6_growth} appears in connection to the weak-type $(1,1)$ estimate for the centered Hardy-Littlewoord maximal operator. 

Another consequence of \zcref{Thm4} is the following positive result. 

\begin{corollary}\label{Cor6}
\zcref[S]{Conj3} is true for non-negative and radial decreasing functions $f$ with $C_d \sim \sqrt{d} \, \log d$.
\end{corollary} 

Our final result is the counterpart of \zcref{thm2} for radial functions. 

\begin{theorem}\label{Thm7}
\zcref[S]{Conj3} is true for radial characteristic functions $f$ with $C_d \sim \sqrt{d} \, \log d$.
\end{theorem}

We highlight the fact that \zcref{Thm7} not only establishes the validity of \zcref{Conj3}, but it does so with the sharp possible dependency of the constant on the dimension $d$, in light of \zcref{Cor1.6_growth}. 

As discussed at the end of \zcref{SubSec1.1}, in the setup of \zcref{Conj3}, an intermediate question would be to establish the validity of \zcref{20260213_17:41} with a constant $C_{d, \varphi}$ on the right-hand side. \zcref[S]{Cor6,Thm7} provide evidence that such a constant should be uniform on $\varphi$. It is also a plausible hypothesis that $C_d$ could be taken to be the optimal constant from case of the Hardy-Littlewood maximal operator.

\subsection{Notation} Vectors in $\R^d$ are denoted here with bold font (e.g., $\x$, $\y$, $\z$) and numbers in $\R$ with regular font (e.g., $x,y,z$). If $E \subset \R^d$, its characteristic function is denoted by $\chi_E$ and its Lebesgue measure by $|E|$. We write $A \lesssim B$ if $A \leq C B$ for a certain constant $C>0 $, and we write $A \sim B$ if $A \lesssim B$ and $B \lesssim A$ (parameters of dependence of such a constant $C>0$ appear as a subscript in the inequality sign).
We sometimes denote a ball centered in $\vt 0$ by $B_r=B(\vt 0,r)$, let $\omega_d$ denote the $d$-dimensional volume of the unit ball $B_1 \subset \R^d$ and let $\sigma_{d-1} = d\, \omega_d$ denote the surface area of the unit sphere $\partial B_1 \subset \R^d$. 

\subsection{AI statement}
No AI was used in the preparation of this paper. 

\subsection{Acknowledgements}

J.W.\ was funded by the European Union's Horizon research and innovation programme under the Marie Skłodowska-Curie Action 2023 No 101151034 (SRMF).

\section{Proof of \texorpdfstring{\zcref{thm2}}{Theorem~\ref{thm2}}}
\subsection{Preliminaries} We may assume without loss of generality that our characteristic function $f$ is such that $f(x) = \liminf_{y\to x}f(y)$, see \zcref{representativeassuminglowerlimit}.
Hence, there are numbers $-\infty \leq a_1 < b_1 < a_2 < b_2 < \ldots < a_N < b_N \leq \infty$ such that 
$$f(x)= \sum_{n=1}^{N} \chi_{(a_n, b_n)}(x).$$
Note that $\M_{\varphi}f(x) = f(x) =1$ for $x \in \cup_{n=1}^N (a_n, b_n)$. If $-\infty < a_1$, note that $\M_{\varphi}f$ is increasing in $(-\infty, a_1+\varepsilon]$, for $\varepsilon>0$ small, and hence 
$$\var_{(-\infty, a_{1} + \varepsilon]}(\M_{\varphi}f) = 1 = \var_{(-\infty, a_{1}+\varepsilon]}(f).$$ 
Similarly, if $b_N <\infty$, $\M_{\varphi}f$ is decreasing in $[b_N - \varepsilon, \infty)$, for $\varepsilon>0$ small, and 
$$\var_{[b_N - \varepsilon, \infty)}(\M_{\varphi}f) = 1 = \var_{[b_N- \varepsilon, \infty)}(f).$$ 
The case when $f$ consists of just one interval follows immediately. In the general case, the conclusion of the proof relies in the following local bound between intervals (for $1 \leq n \leq N-1$)
\begin{equation}\label{20260216_15:22}
\var_{[b_n -\varepsilon , a_{n+1}+ \varepsilon]}(\M_{\varphi}f) \leq 2 = \var_{[b_n- \varepsilon, a_{n+1}+\varepsilon]}(f),
\end{equation}
for $\varepsilon >0 $ small.

\subsection{The local bound}

Given the translation and dilation symmetries, in order to prove \zcref{20260216_15:22} we may assume that $b_n = -1$ and $a_{n+1} = 1$.
Again, by symmetry, it suffices to prove
\[
\var_{[0, 1+\varepsilon]}(\M_{\varphi}f) \leq 1
,
\]
which, since
\[
\var_{[1, 1+\varepsilon]}(\M_{\varphi}f)
=
1-\M_\varphi f(1)
,
\]
is equivalent to
\begin{equation} \label{20260508_15:47}
\var_{[0, 1]}(\M_{\varphi}f) \leq \M_\varphi f(1)
,
\end{equation}
provided $f=0$ on $[-1,1]$.

\subsubsection{Proof by finite differences}
\label{Sec2.2.1_derCon}

The following lemma is the heart of the matter.
\begin{lemma}\label{Lemma_Reg_dim_1}
Let $-1\leq x<y\leq 1$.
Then
\begin{align}\label{20260417_17:32}
\big| \M_\varphi f(y) - \M_\varphi f(x)\big| \leq y\M_\varphi f(y) - x\M_\varphi f(x).
\end{align}
\end{lemma}

\begin{proof}
We first prove the upper bound in \zcref{20260417_17:32}. It clearly holds when $y=1$, so we may assume that $-1 \leq x < y < 1$.
Let $r>0$.
The plan is to compare $f*\varphi_r(y)$ with $f*\varphi_{\tilde r}(x)$ for
\(
\tilde r
=
\big(\frac{1-x}{1-y}\big)r
>
r.
\)
For all $t \leq -1$ we have
\begin{align*}
0\leq   \frac{x-t}{\tilde r} \leq \frac{x-t}r \leq \frac{y-t}r.
\end{align*}
For all $t\geq1$ we have
\begin{align*}
0 \geq \frac{x-t}{\tilde r} = \left(\frac{x-t}r\right)
\left(\frac{1-y}{1-x}\right) \geq  \left(\frac{x-t}r\right)
\left(\frac{t-y}{t-x}\right) = \frac{y-t}r.
\end{align*}

\begin{figure}[htbp]
\centering
\input{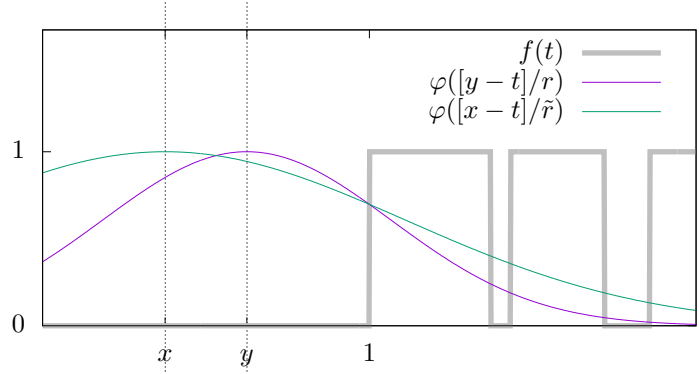}
\caption{The choice of $\tilde r$ satisfies $\varphi\left(\frac{x-t}{\tilde r}\right)\geq \varphi\left(\frac{y-t}{r}\right)$ for all $t\geq1$.}
\end{figure}

Since $f\equiv0$ on $(-1,1)$ and $\varphi$ is decreasing away from $0$, for all $t\in\R$ we can conclude that 
\[
f(t)\,
\varphi\!\left(
\frac{x-t}{\tilde r}
\right)
\geq
f(t)\,
\varphi\!\left(
\frac{y-t}r
\right)
.
\]
Therefore,
\[
\M_\varphi f(x)
\geq
\frac1{\tilde r}
\int_{\R}
f(t)\,
\varphi\!\left(
\frac{x-t}{\tilde r}
\right)
\d t
\geq
\left(\frac{1-y}{1-x}\right)
\frac1r
\int_{\R}
f(t)\,
\varphi\!\left(
\frac{y-t}r
\right)
\d t.
\]
Taking the supremum over all $r>0$ we conclude that 
\[
\M_\varphi f(x)
\geq
\left(\frac{1-y}{1-x}\right)
\M_\varphi f(y)\,,
\]
which is equivalent to the desired upper bound
\[
\M_\varphi f(y) - \M_\varphi f(x)
\leq
y\M_\varphi f(y) - x\M_\varphi  f(x).
\]
To conclude the lower bound define $\tilde{f} = f(-\, \cdot)$, $\tilde{\varphi} = \varphi(-\, \cdot)$, $\tilde x=-y$ and $\tilde y=-x$.
Then \(\M_{\tilde{\varphi}}(\tilde{f})(\tilde{x})=\M_\varphi f(-\tilde{x})=\M_\varphi f( y)\), and we obtain the lower bound in \zcref{20260417_17:32} by inserting those quantities for $f,\varphi,x,y$ in the upper bound.
This concludes the proof of \zcref{20260417_17:32}.
\end{proof}

For a partition $0=x_0<x_1<\ldots<x_n =1$, we use \zcref{20260417_17:32} to get a telescoping sum
\[
\sum_{k=1}^n
\big|
\M_\varphi f(x_k)
-
\M_\varphi f(x_{k-1})
\big|
\leq
\sum_{k=1}^n
\big(x_k \M_\varphi f(x_k)
-
x_{k-1}\M_\varphi f(x_{k-1})\big) = \M_\varphi f(1)
,
\]
which yields the desired bound \zcref{20260508_15:47}.

\subsubsection{Proof by derivatives}

We next present an alternative, possibly more instructive proof of \zcref{20260508_15:47}.
The starting point is the following corollary of \zcref{Lemma_Reg_dim_1}:

\begin{figure}[htbp]
\centering
\input{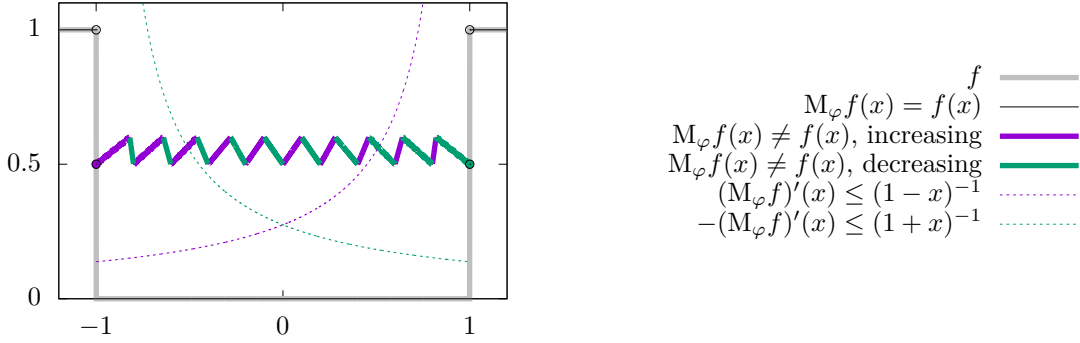}
\caption{Possible behaviour of the maximal function in the disconnecting set.}
\end{figure}

\begin{corollary}\label{Cor_Lipschitz}
The function $\M_{\varphi}f$ is locally Lipschitz in $(-1,1)$ and for almost all $ x \in (-1,1)$ we have the pointwise bound
\begin{align*}
-\frac{\M_{\varphi}f (x)}{1+x} \leq \big(\M_{\varphi}f\big)'(x) \leq \frac{\M_{\varphi}f (x)}{1-x}.
\end{align*}
\end{corollary}

\begin{proof}
\zcref[S]{20260417_17:32} can be reformulated as
\[
-\frac{\M_{\varphi}f (x)}{1+y}
\leq
\frac{
\M_\varphi f(y)
-
\M_\varphi f(x)
}{y-x}
\leq
\frac{\M_{\varphi}f (y)}{1-x}\,,
\]
and we finish the proof by letting $y\rightarrow x$.
\end{proof}

We will insert this bound into the following \zcref*[noref,nocap]{Lemma12_Abs_cont}:

\begin{lemma}\label{Lemma12_Abs_cont}
Let $g:[0,1) \to \R_{\geq 0}$ be a locally absolutely continuous function such that 
\begin{equation}\label{20260508_17:06}
\frac{-g(x)}{1+x} \leq g'(x) \leq \frac{g(x)}{1-x}
\end{equation}
for each $x \in [0,1)$. Then, for each $x \in [0,1)$, 
$$\int_0^x |g'(y)|  \d y \leq x \, g(x).$$
\end{lemma}
\begin{proof}
For $x \in [0,1)$ let 
$$F(x)= x \, g(x) - \int_0^x |g'(y)| \d y.$$
The lemma is equivalent to showing that $F(x) \geq 0$. Note that $F(0) = 0$. We claim that $F'(x) \geq 0$ for all $x \in (0,1)$. In fact, by the fundamental theorem of calculus: 
$$F'(x) = g(x) + x g'(x) - |g'(x)|.$$
If $g'(x) \geq 0$, we have $F'(x) \geq 0$ if and only if $ g'(x)\leq  \frac{g(x)}{1-x}$. If $g'(x) \leq 0$, we have $F'(x) \geq 0$ if and only if $g'(x)\geq  \frac{-g(x)}{1+x}$, and these are exactly the conditions in \zcref{20260508_17:06}.
\end{proof}
From \zcref{Cor_Lipschitz} and \zcref{Lemma12_Abs_cont} we get, for $0 < x < 1$, 
$$\var_{[0, x]}(\M_{\varphi}f)  = \int_0^x \big|(\M_{\varphi}f)'(y)\big|  \d y \leq x \M_\varphi f(x).
$$
Sending $x \to 1^-$ and using \zcref{Prop11_cont} (or \zcref{Prop11_contbv}) we get \zcref{20260508_15:47}.

\begin{remark}
One can deduce \(\M_\varphi f(1)\geq\limsup_{x\rightarrow1^-}\M_\varphi f(x)\) from \zcref{Lemma_Reg_dim_1} which would also suffice to complete the proof instead of invoking \zcref{Prop11_cont} or \zcref{Prop11_contbv}.
\end{remark}

\subsection{A remark on the kernel} Note that the proof in this section does not require that $\varphi$ is an even function; it just uses the fact that $\varphi$ is decreasing away from $0$. Our \zcref{thm2} holds for this slightly larger class of convolution kernels in dimension $d=1$.

\section{Proof of \texorpdfstring{\zcref{Thm7}}{Theorem \ref{Thm7}}}

\subsection{Setup} We may assume here without loss of generality that $d \geq 2$ since the case $d=1$ was already treated in \zcref{thm2}. Let $f\in\BV(\mathbb{R}^d)$ be a radial characteristic function, i.e., there exists $\rad f:[0,\infty)\rightarrow\mathbb{R}$ with \(f(\vt x)=\rad f(|\vt x|)\). Then, by \zcref{varfradial}, $\rad f\in\BV_\loc((0,\infty))$; more precisely, for each $r>0$ we have \(\rad f \cdot \chi_{(r,\infty)} \in\BV((0,\infty))\). Moreover, it suffices to assume $f(\vt 0)=0$ and, by \zcref{radialcharfrepresentative}, that there exist $N\in\N\cup\{\infty\}$ and sequences $\{a_n\}_{-N-1<n\leq0}$ and $\{b_n\}_{-N-1<n\leq0}$, with $0< a_n<b_n<a_{n+1}$ for each $n$, such that 
\[
\rad f = \sum_{n=-N}^0 \chi_{(a_n, b_n)}.\]
In this case we observe that 
\begin{align}\label{20260826_08:36}
\var f
=
\sigma_{d-1}\sum_{n=-N}^0 \left(a_n^{d-1}+b_n^{d-1}\right)\,,
\end{align}
and we are assuming that this is a finite quantity. Note that the reason we are indexing the sequence of intervals with negative integers is that, even with \zcref{20260826_08:36} finite, it may contain an infinite number of intervals that accumulate near $0$.

If $N = \infty$, note that $\inf_{n\leq 0}\{a_n\} = 0$. If $N <\infty$ we set $b_{-N-1} :=0$. Since $\M_\varphi f$ is also radial, there exists $\rad{\M_\varphi f}:[0,\infty)\rightarrow\R$ with $\M_\varphi f(\vt x)=\rad{\M_\varphi f}(|\vt x|)$. Noting that $\rad{\M_\varphi f}$ is constant on $(a_n,b_n)$, from \zcref{approx_by_annuli,varfradial} we have
\begin{align}\label{20260903_14:36}
\begin{split}
\var \M_\varphi f
&=
\sigma_{d-1}
\int_0^\infty
r^{d-1}
\d\left|\grad\rad{\M_\varphi f}\right|(r)
\\
&=
\sum_{n=-N}^0
\sigma_{d-1} \int_{[b_{n-1},a_n]}
r^{d-1}
\d\left|\grad\rad{\M_\varphi f}\right|(r)
+
\sigma_{d-1} \int_{[b_0,\infty)}
r^{d-1}
\d\left|\grad\rad{\M_\varphi f}\right|(r)
\\
&\leq
\sum_{n=-N}^0
\sigma_{d-1} \, a_n^{d-1}
|\grad\rad{\M_\varphi f}|([b_{n-1},a_n])
\ +
\sigma_{d-1} \int_{[b_0,\infty)}
r^{d-1}
\d\left|\grad\rad{\M_\varphi f}\right|(r)\,,
\end{split}
\end{align}
where $\grad\rad{\M_\varphi f}$ denotes the distributional derivative of $\rad{\M_\varphi f}$ and $\big|\grad\rad{\M_\varphi f}\big|$ denotes its total variation measure. 

\begin{remark}
\label{bound_classical_instead_of_weak}
In \zcref{20260903_14:36} the term
\(
|\grad\rad{\M_\varphi f}|([b_{n-1},a_n])
=
\var_{(b_{n-1}-\varepsilon_n,a_n+\varepsilon_n)}\rad{\M_\varphi f}
\)
denotes the variation in terms of its weak definition.
We will instead bound the classical variation of $\rad{\M_\varphi f}$, as defined in \zcref{definition1dvariation} via finite differences, and use that it bounds the weak definition used in \zcref{20260903_14:36}, see \zcref{appendix_bvoned}.
\end{remark}

\subsection{Derivative control}

Here we establish the key \zcref{Lemma3.1_Em} used to control the term
\(
\var_{(b_{n-1}-\varepsilon_n,a_n+\varepsilon_n)}\rad{\M_\varphi f}
.
\)
To that end we shall need estimates in the spirit of the ones in \zcref{Sec2.2.1_derCon} to control the derivative of the maximal function. In order to minimize the technicalities we argue here via finite differences.

Let us first establish some notation. Abbreviate $\vt e=(1,0,\ldots,0)\in\mathbb{R}^d$.
For $g:\mathbb{R}^d\rightarrow\mathbb{R}$ and $\delta\in\mathbb{R}$ denote by $\Delta^\delta g:\mathbb{R}^d\rightarrow\mathbb{R}$ the map
\[
\Delta^\delta g(\vt x)
=
g(\vt x+\delta\vt e)
-
g(\vt x)
.
\]
Then, recalling
\[
g*h(\vt x)
=
\int_{\mathbb{R}^d}
g(\vt x-\vt y)
h(\vt y)
\intd\vt y
=
\int_{\mathbb{R}^d}
g(\vt y)
h(\vt x-\vt y)
\intd\vt y
,
\]
we have
\[
\Delta^\delta(g*h)(\vt x)
=
(\Delta^\delta g)*h(\vt x)
=
g*(\Delta^{\delta} h)(\vt x)
\]
and
\begin{align}\label{20260905_19:10}
\Delta^{\delta}\chi_{B(\vt 0,r)}(\vt x-\vt y)
=
\chi_{B(\vt x+\delta\vt e,r)}(\vt y)
-
\chi_{B(\vt x,r)}(\vt y)
.
\end{align}

Define the hyperplane 
\[
H(\delta)
:=
\{\vt x \in \R^d \ :\ 
\vt x\cdot\vt e\geq\delta
\}
.
\]

\begin{figure}[htbp]
\centering
\includegraphics{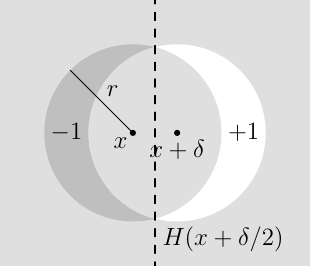}
\caption{The regions of negative and positive signs of $\vt{y}\mapsto \Delta^{\delta} \chi_{B(\vt 0,r)}(\x - \y)$.}
\label{ballsign}
\end{figure}  

A simple yet useful observation for our purposes is the following: for $\delta > 0$ and $r >0$ we have \( \Delta^\delta \chi_{B(\vt 0,r)}(x\vt e-\vt y)\geq0\) for $\vt y\in H(x+\delta/2)$ and $ \Delta^\delta \chi_{B(\vt 0,r)}(x\vt e-\vt y)\
\leq0$ for $\vt y\in\mathbb{R}^d\setminus H(x+\delta/2)$; see \zcref{ballsign}.

\begin{figure}[htbp]
\begin{subfigure}{0.5\textwidth}
\centering
\includegraphics{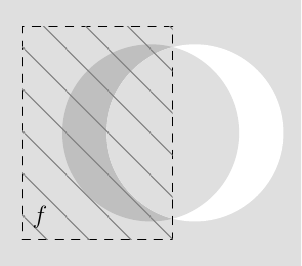}
\caption{Shifting the ball.}
\end{subfigure}
\begin{subfigure}{0.49\textwidth}
\centering
\includegraphics{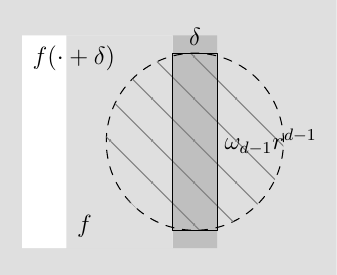}
\caption{Shifting the rectangle $f$.}
\end{subfigure}
\caption{A function $f$ for which $\Delta^\delta(f*\chi_{B(\vt 0,r)})(x\vt e)$ is maximal.}
\label{largestchange}
\end{figure}

Thus, for any characteristic function $f$, and any $\delta > 0$ and  $x\geq0$, we have
\begin{align*}
f*(\Delta^\delta\chi_{B(\vt 0,r)})(x\vt e)
&\leq \chi_{H(x+\delta/2)}*(\Delta^\delta\chi_{B(\vt 0,r)})(x\vt e)
= (\Delta^\delta \chi_{H(x+\delta/2)})*\chi_{B(\vt 0,r)}(x\vt e) \leq \delta\,\omega_{d-1}r^{d-1},
\end{align*}
see \zcref{largestchange}.
By symmetry we can conclude that for all $\delta >0 $ and $x \geq 0$
\begin{equation}
\label{gradientawayfromsupport1}
\left|
\Delta^\delta\big(f*\chi_{B(\vt 0,r)}/(\omega_d\,r^d)\big)(x\vt e)
\right|
\leq
\frac{\delta\, \omega_{d-1}}{\omega_d\, r}
\sim
\frac{\delta
\sqrt d
}r
.
\end{equation}

For a generic function $g$, not necessarily non-negative, let
\begin{equation}\label{20261009_03:10}
\M_{\varphi}^{s}g(\x) := \sup_{t>0}\big(g*\varphi_t\big)(\x)
\end{equation}
denote the maximal operator of convolution type without the absolute value inside (the superscript s here would stand for `signed'), and set $\Ms$ for the analogous signed version of the Hardy-Littlewood maximal operator.  For our type of kernels $\varphi$ (non-negative, radial decreasing and with integral $1$), by the layer cake representation we still have
\begin{align}\label{20260903_10:54}
\Ms_{\varphi}g(\x) \leq \Ms g(\x).
\end{align}
as in \zcref{20260903_10:40}. We now move on to the key lemma. 
\begin{lemma}\label{Lemma3.1_Em}
With notations as above, if $b_{n-1} \leq x < a_n$ and $0 < \delta \leq  a_n - x$ we have 
\begin{align}\label{20260903_11:10}
- \,\frac{\delta\sqrt d}{x - b_{n-1}} -  \,\frac{\delta\sqrt d}{ \sqrt{a_n^2 - \big(x + \frac{\delta}{2}\big)^2}}  \ \lesssim \ \Delta^\delta(\M_{\varphi} f)(x\vt e) \  \lesssim \ \frac{\delta\sqrt d}{a_n-x - \delta}.
\end{align}
\end{lemma}
\begin{proof} By the sublinearity of the maximal operators and \zcref{20260903_10:54}, we have
\begin{align}\label{20260903_11:09}
\begin{split}
- \Ms \left(-\Delta^\delta f\right)(x\vt e)   \leq - \Ms_{\varphi} \left(-\Delta^\delta f\right)(x\vt e)  & \leq \Delta^\delta(\M_{\varphi} f)(x\vt e)  =  \Delta^\delta(\Ms_{\varphi} f)(x\vt e) \\
 & \leq \Ms_{\varphi} \left(\Delta^\delta f\right)(x\vt e) \leq   \Ms \left(\Delta^\delta f\right)(x\vt e).
\end{split}
\end{align}
Hence, in light of \zcref{20260903_11:09}, it suffices to consider the signed Hardy-Littlewood maximal operator in order to establish \zcref{20260903_11:10}.

Write \(f= f_1 + f_2\), where \(f_1 = f\cdot \chi_{B(\vt 0,b_{n-1})}\) and \(f_2 = f\cdot \chi_{\mathbb{R}^d\setminus B(\vt 0,a_n)}\).
\begin{figure}[htbp]
\centering
\includegraphics{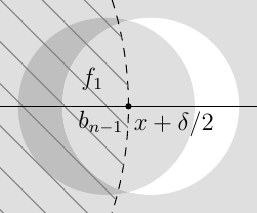}
\caption{For $f_1 = f\cdot \chi_{B(\vt 0,b_{n-1})}$ and $x+\delta/2\geq b_{n-1}$ we have $f_1*(\Delta^{\delta}\chi_{B(\vt 0,r)})(x\vt e)\leq0$.}
\label{decrease}
\end{figure}
Recalling \zcref{20260905_19:10}, note that for any $r >0$
\begin{align}\label{20260903_11:23}
\begin{split}
\big(\Delta^\delta f_1*\chi_{B(\vt 0,r)}\big)(x\vt e)
& = f_1*\big(\Delta^\delta \chi_{B(\vt 0,r)}\big)(x\vt e) \\
& = \int_{\R^{d}} f_1(\y) \left( \chi_{B((x+\delta) \vt e, r)}(\y) -  \chi_{B(x \vt e, r)}(\y)\right) \d \y \leq 0\, ,
\end{split}
\end{align}
see \zcref{decrease}.
Also, note that $\big(\Delta^{\delta}f_2*\chi_{B(\vt 0,r)}\big)(x\vt e) = 0$ if $r < a_n - x - \delta$. Hence, in order to get a non-trivial average of $\Delta^\delta f_2$ one must consider a radius $r \geq a_n - x - \delta$. From \zcref{20260903_11:23} and \zcref{gradientawayfromsupport1} we get 
\begin{align*}
 \Ms \left(\Delta^\delta f\right)(x\vt e) \leq \Ms \left(\Delta^\delta f_2\right)(x\vt e)\lesssim \frac{\delta\sqrt d}{a_n-x - \delta}.
\end{align*}
This establishes the inequality on the right-hand side of \zcref{20260903_11:10}.

Similarly, observe that $\big(\Delta^\delta f_1*\chi_{B(\vt 0,r)}\big)(x\vt e) = 0$ if $r <  x - b_{n-1}$. From \zcref{gradientawayfromsupport1} we plainly get 
\begin{align}\label{20260903_12:42}
 -\Ms \left(-\Delta^\delta f_1\right)(x\vt e) \gtrsim - \,\frac{\delta\sqrt d}{x - b_{n-1}}.
\end{align}
\begin{figure}[htbp]
\centering
\includegraphics{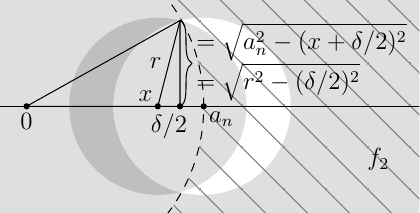}
\caption{For \(f_2 = f\cdot \chi_{\mathbb{R}^d\setminus B(\vt 0,a_n)}\), if $r\leq\sqrt{(\delta/2)^2+a_n^2-(x+\delta/2)^2}$ then $f_2*(\Delta^{\delta}\chi_{B(\vt 0,r)})(x\vt e)\geq0$.}
\label{sqrtincrease}
\end{figure}

\noindent Finally, observe that if $r < \sqrt{\big(\frac{\delta}{2}\big)^2+ a_n^2 - \big(x + \frac{\delta}{2}\big)^2}$, then $B(x\vt e, r) \setminus H\big(x + \frac{\delta}{2}\big)$ does not intersect $\R^d \setminus B(\vt 0,a_n)$ and hence, recalling \zcref{20260905_19:10},
\begin{align*}
\big(-\Delta^\delta f_2*\chi_{B(\vt 0,r)}\big)(x\vt e) & = f_2*\big(-\Delta^\delta \chi_{B(\vt 0,r)}\big)(x\vt e)\\
& =  \int_{\R^{d}} f_2(\y) \left( \chi_{B(x \vt e, r)}(\y) -  \chi_{B((x+ \delta) \vt e, r)}(\y)\right) \d \y \leq 0\,,
\end{align*}
see \zcref{sqrtincrease}.
From \zcref{gradientawayfromsupport1} we then get 
\begin{align}\label{20260903_12:43}
 -\Ms \left(-\Delta^\delta f_2\right)(x\vt e) \gtrsim - \,\frac{\delta\sqrt d}{ \sqrt{\big(\frac{\delta}{2}\big)^2+ a_n^2 - \big(x + \frac{\delta}{2}\big)^2}} \geq - \,\frac{\delta\sqrt d}{ \sqrt{ a_n^2 - \big(x + \frac{\delta}{2}\big)^2}}.
\end{align}
Putting together \zcref{20260903_12:42} and \zcref{20260903_12:43} we arrive at the inequality on the left-hand side of \zcref{20260903_11:10}. 
\end{proof}

\begin{remark}
It follows from \zcref{Lemma3.1_Em} that $\M_{\varphi} f$ is locally Lipschitz on each annulus $\{b_{n-1} < |\x| < a_n\}$. 
\end{remark}
 
\subsection{Conclusion of the proof} 

\subsubsection{The components between annuli}
Looking at \zcref{20260903_14:36,bound_classical_instead_of_weak}, we first bound  $\var_{(b_{n-1}-\varepsilon_n,a_n+\varepsilon_n)}\rad{\M_\varphi f}$ for each fixed $n<0$.
First, note
\begin{equation}
\label{initial_and_final_jump}
\var_{(b_{n-1}-\varepsilon_n,b_n]}\M_\varphi f
+
\var_{[a_n,a_n+\varepsilon_n)}\M_\varphi f
=
1-\M_\varphi f(b_{n-1})
+
1-\M_\varphi f(a_n)
\leq2
.
\end{equation}
We proceed to bound $\var_{[b_{n-1},a_n]}\rad{\M_\varphi f}$.
Denote $c_n=(a_n+b_{n-1})/2$ the midpoint of the interval. 

Let us first bound $\var_{[b_{n-1},c_n]}\rad{\M_\varphi f}$. For this, consider a partition $P = \{ b_{n-1} = x_0 < x_1 < \ldots < x_K = c_n\}$ and abbreviate $\Delta_k := \rad{\M_{\varphi} f}(x_{k+1}) - \rad{\M_{\varphi} f}(x_k)$. We can decompose these differences into their positive and negative parts, i.e., $\Delta_k= \big( \Delta_k\big)_+ - \big( \Delta_k\big)_-$. Note that 
$$ -1 \leq  \rad{\M_{\varphi} f}(c_n) - \rad{\M_{\varphi} f}(b_{n-1}) = \sum_{k} \Delta_k =  \sum_{k} \big( \Delta_k\big)_+ -  \sum_{k} \big( \Delta_k\big)_-$$
and therefore 
\begin{align}\label{20260903_15:25}
\sum_{k} \big|\Delta_k\big|  = \sum_{k} \big( \Delta_k\big)_+ +  \sum_k\big( \Delta_k\big)_-  \leq  1 + 2\sum_{k} \big( \Delta_k\big)_+.
\end{align}
From \zcref{Lemma3.1_Em} we have, for $k = 0, 1, \ldots, K-1$,  
$$\big( \Delta_k\big)_+ \lesssim  \frac{(x_{k+1} - x_k)\sqrt d}{a_n- x_{k+1}} \leq \frac{(x_{k+1} - x_k)\sqrt d}{a_n-c_n}.$$
Summing over $k$ and using \zcref{20260903_15:25} we arrive at 
$$\sum_{k} \big|\Delta_k\big|   \lesssim \sqrt{d}\,,$$
and hence
 \begin{align}\label{20260904_17:45}
 \var_{[b_{n-1},c_n]}\rad{\M_\varphi f}\lesssim \sqrt{d}.
 \end{align}
 
Let us now bound $\var_{[c_n, a_n]}\rad{\M_\varphi f}$. Consider then a partition $P = \{ c_n = x_0 < x_1 < \ldots < x_K = a_n\}$. Similarly as before, from the fact that 
$$ 1 \geq  \rad{\M_{\varphi} f}(a_n) - \rad{\M_{\varphi} f}(c_{n}) = \sum_{k} \Delta_k =  \sum_{k} \big( \Delta_k\big)_+ -  \sum_{k} \big( \Delta_k\big)_-$$
one arrives at 
\begin{align}\label{202600908_17:31}
\sum_{k} \big|\Delta_k\big|  = \sum_{k} \big( \Delta_k\big)_+ +  \sum_k\big( \Delta_k\big)_-  \leq  1 + 2\sum_{k} \big( \Delta_k\big)_-.
\end{align}
From \zcref{Lemma3.1_Em} we have, for $k = 0, 1, \ldots, K-1$,  
\begin{align*}
\big( \Delta_k\big)_- & \lesssim  \frac{(x_{k+1} - x_k)\sqrt d}{c_n- b_{n-1}}  + \frac{(x_{k+1} - x_k)\sqrt d}{\sqrt{a_n^2 - \left(\frac{x_k + x_{k+1}}{2}\right)^2}}.
\end{align*}
Summing over $k$ and using \zcref{202600908_17:31} we arrive at 
\begin{align*}
\sum_{k} \big|\Delta_k\big| \lesssim \sqrt{d} + \sqrt d \, \sum_{k}  \frac{(x_{k+1} - x_k) }{\sqrt{a_n^2 - \left(\frac{x_k + x_{k+1}}{2}\right)^2}}.
\end{align*}
When we consider the limit as $|P| := \max\{x_{k+1} -x_k\} \to 0$, this Riemann sum converges to 
\begin{align*}
\int_{c_n}^{a_n} \frac{1 }{\sqrt{a_n^2 - t^2}} \d t  \leq \int_{0}^{a_n} \frac{1 }{\sqrt{a_n^2 - t^2}} \d t = \int_{0}^{1} \frac{1 }{\sqrt{1 - t^2}} \d t  = \frac{\pi}{2}.
\end{align*}
This shows that 
\begin{align}\label{20260904_17:46}
\var_{[c_n, a_n]}\rad{\M_\varphi f} \lesssim \sqrt{d}.
\end{align}
From \zcref{initial_and_final_jump,20260904_17:45,20260904_17:46} we conclude that 
\begin{align}\label{20260904_22:20}
\var_{(b_{n-1}-\varepsilon_n, a_n+\varepsilon_n)}\rad{\M_\varphi f} \lesssim \sqrt{d}.
\end{align}

\subsubsection{The component near infinity} 

We start by recalling that for each $\x \in \R^d$
\begin{equation*}
\M_{\varphi}f(\x) \leq  \M f(\x) \leq \M \chi_{B(\vt 0, b_0)}(\x).
\end{equation*} 
Let $\varepsilon>0$ be small such that $a_0 < b_0 - \varepsilon$, and let $\Omega = \{ \x \in \R^d \ : \ |\x| > b_0 - \varepsilon\}$. From \zcref{varfradial} and the coarea formula (\zcref{coarea}) we have
\begin{align}\label{20260904_22:21}
\begin{split}
\sigma_{d-1} \int_{[b_0,\infty)}
r^{d-1}
\d\left|\grad\rad{\M_\varphi f}\right|(r) & = \var_{\Omega}(\M_{\varphi}f)
=
\int_{-\infty}^\infty
P(\{x\in\mathbb{R}^d:\M_{\varphi}f(x)>t\},\Omega)\, \d t \\
& \leq \int_{-\infty}^\infty
P(\{x\in\mathbb{R}^d:\M f(x)>t\},\Omega)\, \d t\\
& \leq \int_{-\infty}^\infty
P(\{x\in\mathbb{R}^d:\M \chi_{B(\vt 0, b_0)}(x)>t\},\Omega)\, \d t\\
& =  \var_{\R^d}(\M \chi_{B(\vt 0, b_0)})\\
& \sim \sqrt{d} \log d  \, \var_{\R^d}(\chi_{B(\vt 0, b_0)})\\
& = \sqrt{d} \log d \ \sigma_{d-1} b_0^{d-1}\,,
\end{split}
\end{align}
where we have used \zcref{Cor5} in the penultimate passage.

\subsubsection{Conclusion} 
Recalling \zcref{20260826_08:36,bound_classical_instead_of_weak} and plugging the estimates \zcref{20260904_22:20} and \zcref{20260904_22:21} into \zcref{20260903_14:36} we arrive at the desired result. 

\section{Asymptotics for the maximal function of the unit ball}
\label{unitballsection}

In this \zcref*[noref,nocap]{unitballsection} we prove \zcref{Thm4}. Our proof requires a certain elaboration, and hence we divide it into a few components in order to ease the exposition.

\subsection{The volume of a cap}\label{Sec4.1_cap}

For $t\in[-1,1]$ denote by $C(t)$ the volume of the cap $\{\y = (y_1, y_2, \ldots, y_d) \in B_1\ :\ y_1\geq t\}$ divided by $\omega_d=|B_1|$, and let $g:[-1,1] \to \R$ be given by 
\[
g(s)
:=
(1-s^2)^{\frac12}
.
\]
Then
\begin{equation}\label{20260805_14:52}
C(t)
=
\frac{\omega_{d-1}}{\omega_d}
\int_t^1
g(s)^{d-1}
\d s
.
\end{equation}

Our starting point is an algebraic expression that approximates $C(t)$ with an explicit dependence on the dimension $d$.

\begin{lemma}
\label{pro_capvolume}
For $t\geq0$ we have
\[
C(t)
\sim
\frac{
(1-t^2)^{\frac{d+1}2}
}{
t \sqrt d +1
}
,
\]
with the implicit constant being universal. 
\end{lemma}

\begin{proof}The result is trivial for $d=1$, so we may assume that $d>1$ here. Observe that the integrand $g(s)^{d-1}$ in \zcref{20260805_14:52} is a decreasing function for $s\geq 0$.
The idea is to show that the whole integral in $C(t)$ is dominated by the integral from $t$ to $t+\varepsilon(t)$, with $\varepsilon(t)$ such that $g(t+\varepsilon(t))^{d-1}$ is approximately half of $g(t)^{d-1}$.
Then we will be able to conclude $C(t)\sim(\omega_{d-1}/\omega_d)\varepsilon(t)g(t)^{d-1}$.
Note that $g^{d-1}$ approximately halves if and only if $g$ shrinks by a factor of $1-1/d$.
For $0 \leq t <1$, in order to find $\varepsilon(t)$ that achieves this, we estimate the derivative of $g$ on $[t,t+\eta]$, where 
\begin{equation}\label{20260610_17:39}
\eta:=\min\left\{d^{-\frac12}\,,\,\frac{1-t}2\right\}
\end{equation}
is an upper bound for $\varepsilon(t)$, as we will see. 
Since
\(
g'(s)
=
-\frac s{(1-s^2)^{1/2}}
\)
is decreasing, we can take $g'(t)$ as an upper bound for $g'$ on $[t,t+\eta]$.
We will show that
\[
h(t)
:=
-\frac{2\big(t+d^{-\frac12}\big)}{(1-t^2)^{\frac12}}
,
\]
which is decreasing as well, is a lower bound for $g'$ on $[t,t+\eta]$.
To that end, it suffices to estimate
\begin{align*}
\frac{g'(t+\eta)}{h(t)}
=
\frac{t+\eta}{2\big(t+d^{-\frac12}\big)}
\left(
\frac{1-t}{1-t-\eta}\right)^{\frac12}\left(\frac{1+t}{1+t+\eta}
\right)^{\frac12}
\leq \frac{t+\eta}{2\big(t+d^{-\frac12}\big)}
\left(
\frac{1-t}{1-t-\eta}\right)^{\frac12} \leq  \frac{\sqrt{2}}{2} < 1\,,
\end{align*}
where we have used \zcref{20260610_17:39}. We can now define
\[
\varepsilon(t)
:=
 - \frac{
g(t)
}{
d\, h(t)
}
=
\frac{
1-t^2
}{
2d\big(t+d^{-\frac12}\big)
}
,
\]
and one can check that 
\(
\varepsilon(t)
\leq \eta
\).
Hence, for all $t\leq s\leq t+\varepsilon(t)$, we conclude that  
\begin{equation}
\label{gdirectlowerbound}
g(s)-g(t)
\geq
(s-t)g'(s)
\geq
\varepsilon(t)h(t)
=
-\frac{
g(t)
}d
.
\end{equation}
From \zcref{gdirectlowerbound} and the elementary inequality 
\begin{equation}\label{20260805_14:14}
\left(1 - \frac{1}{d}\right)^{d-1} \geq e^{-1}\,,
\end{equation} 
valid for $d>1$, we get 
\begin{equation}\label{20260805_14:02}
\frac{\omega_d}{\omega_{d-1}}
C(t)
\geq
\int_t^{t+\varepsilon(t)}
g(t)^{d-1}
\left(
1-\frac1d
\right)^{d-1}
\d s
\geq
e^{-1}\,\varepsilon(t)\,g(t)^{d-1}.
\end{equation}
This plainly yields
\[
C(t)
\gtrsim
\frac{
(1-t^2)^{\frac{d+1}2}
}{
t \sqrt d +1
}
.
\]

To prove the upper bound, note that 
\(
\max\{\varepsilon(t),t\}
\geq
\frac1{
3\sqrt{d}
}\,,
\)
which implies \(t+\frac{\varepsilon(t)}{2}\gtrsim t+d^{-\frac12}\). Hence, there exists a universal $c$, with $0<c < 1$, such that 
$$-g'\Big(t+\tfrac{\varepsilon(t)}{2}\Big)\geq-2\,c\,h(t).$$
From this we get
\[
g(t)-g\big(t+\varepsilon(t)\big)
\geq
g\Big(t+\tfrac{\varepsilon(t)}{2}\Big)-g\big(t+\varepsilon(t)\big)
\geq
-
\frac{
\varepsilon(t)
}2
g'\Big(t+\tfrac{\varepsilon(t)}{2}\Big)
\geq
-c\,\varepsilon(t)
\,h(t)
=
\frac{
c\,g(t)
}d\,,
\]
which means
\begin{equation}\label{20260819_21:34}
g\big(t+\varepsilon(t)\big)
\leq
\left(1 - \frac{c}{d} \right)
g(t)
.
\end{equation}
Since $\varepsilon(t)$ and $g(t)$ are decreasing in $t$, we can conclude inductively that 
\begin{align}\label{20260610_21:39}
g\big(t+k\varepsilon(t)\big)
\leq
\left(
1
-
\frac cd
\right)
^k
g(t)
\end{align}
for all $k \geq 0$ such that $t+k\varepsilon(t) \leq 1$. Using again the fact that $g$ is decreasing, splitting the integral defining $C(t)$ into blocks of size $\varepsilon(t)$, and invoking \zcref{20260610_21:39}, we arrive at
\begin{align}\label{20260610_21:34}
\frac{\omega_d}{\omega_{d-1}}
C(t)
&\leq
\sum_{k=0}^{\lfloor(1-t)/\varepsilon(t)\rfloor }
\varepsilon(t)
\,g\big(t+k\varepsilon(t)\big)^{d-1}
\leq
\varepsilon(t)\, g(t)^{d-1}
\sum_{k=0}^{\lfloor(1-t)/\varepsilon(t)\rfloor }
\left(
1
-
\frac cd
\right)
^{(d-1)k}.
\end{align}
From the elementary inequality $(1 - \frac{c}{d})^{d} \leq e^{-c}$ we get 
\begin{align}\label{20260610_21:35}
\sum_{k=0}^\infty
\left(
1
-
\frac cd
\right)
^{(d-1)k}
=
\sum_{k=0}^\infty
\left(
1
-
\frac cd
\right)
^{d\frac{(d-1)k}{d}}
\leq
\sum_{k=0}^\infty
e^{-c\frac{(d-1)k}{d}}
\leq
\sum_{k=0}^\infty
e^{-ck/2}
\lesssim
1
.
\end{align}
From \zcref{20260610_21:34} and \zcref{20260610_21:35} we conclude that 
\[C(t)
\lesssim
\frac{
(1-t^2)^{\frac{d+1}2}
}{
t\sqrt d +1
}
.
\]
\end{proof}

The next result bounds how much our function $C$ can decrease locally. 

\begin{lemma}
\label{lem_capvolumedisturbance}
As before, let 
\begin{equation}\label{20260807_16:53}
\varepsilon(t)
=
\frac{
1-t^2
}{
2d\big(t+d^{-\frac12}\big)
}.
\end{equation}
Let $0\leq t<s<1$, with $s\leq 2t$ and $1-s^2\geq(1-t^2)/2$.
Then
\begin{equation}\label{20260805_14:09}
C(s)
\gtrsim
e^{-\frac{8(s-t)}{\varepsilon(t)}}
C(t)
.
\end{equation}
\end{lemma}

\begin{proof} If $d=1$ the verification is straightforward. We henceforth assume $d>1$. Note that the assumptions on $s$ and $t$ imply that $\varepsilon(t)\geq \varepsilon(s)\geq\varepsilon(t)/4$.

Let $t_0:=t$ and $t_{k+1}:=t_k+\varepsilon(t_k)$.
By induction, \zcref{gdirectlowerbound} and \zcref{20260819_21:34} imply
\begin{equation}\label{20260819_21:38}
\left(1-\frac{c}{d}\right)^k
g(t) \geq g(t_k)
\geq
\left(1-\frac1d\right)^k
g(t)
.
\end{equation}
In particular, note from \zcref{20260819_21:38} that $t_k \to 1$ as $k \to \infty$. Take $k$ such that $t_k\leq s< t_{k+1}$. If $k =0$,
\begin{equation}\label{20260819_21:46}
g(s) \geq g(t_1) \geq \left(1-\frac1d\right)
g(t),
\end{equation}
and from \zcref{20260805_14:02}, \zcref{20260819_21:46}, \zcref{20260610_21:34} and \zcref{20260610_21:35} we get
\begin{align*}
\frac{\omega_d}{\omega_{d-1}}
C(s)
&\gtrsim
g(s)^{d-1}\varepsilon(s)
\gtrsim
g(t)^{d-1}\varepsilon(t)
\gtrsim
\frac{\omega_d}{\omega_{d-1}}
C(t)\,,
\end{align*}
which is \zcref{20260805_14:09} in this case. 

If $k\geq1$, note that 
\begin{equation}\label{20260819_21:53}
s-t\geq t_k-t\geq k \varepsilon(s) \geq k\varepsilon(t)/4\geq(k+1)\varepsilon(t)/8\,,
\end{equation}
and thus, from \zcref{20260819_21:38} and \zcref{20260819_21:53}, 
\begin{align*}
g(s)
&\geq
g(t_{k+1})
\geq
\left(1-\frac1d\right)^{k+1}
g(t)
\geq
\left(1-\frac1d\right)^{\frac{8(s-t)}{\varepsilon(t)}}
g(t)
\geq
e^{-\frac{8(s-t)}{(d-1)\varepsilon(t)}}
g(t)\,,
\end{align*}
where we have used \zcref{20260805_14:14} in the last inequality. The previous estimate, in light of \zcref{20260805_14:02}, \zcref{20260610_21:34} and \zcref{20260610_21:35}, yields
\begin{align*}
\frac{\omega_d}{\omega_{d-1}}
C(s)
&\gtrsim
g(s)^{d-1}\varepsilon(s)
\gtrsim
e^{-\frac{8(s-t)}{\varepsilon(t)}}
g(t)^{d-1}\varepsilon(t)
\gtrsim
e^{-\frac{8(s-t)}{\varepsilon(t)}}
\frac{\omega_d}{\omega_{d-1}}
C(t)\,,
\end{align*}
as desired.
\end{proof}

\subsection{Auxiliary maximal functions}

Our goal is to find an estimate for the volume of the intersection of two $d$-dimensional balls, one of which is the unit ball $B_1 = B({\bm 0},1)$.

Without loss of generality, throughout the rest of this section, we consider a point $\x \in \R^d$ in the $1$-axis, i.e. $\x = (x, 0, 0, \ldots, 0)$, where $x = |\x| > 1$. Let $-1< s< 1$.
Then $\partial B_1 \cap\partial B(\x,x+s)$ is a $(d-2)$-dimensional sphere orthogonal to the $1$-axis.
Let $t\in\mathbb{R}$ be its $1$-component.
Then $t$, $s$ and $x$ are related by the following identity
\begin{align}\label{20260806_16:18}
1-t^2
&=
(x+s)^2
-
(x-t)^2
=
2x(s+t)
+
s^2
-
t^2\,,
\end{align}
which is equivalent to 
\begin{align}\label{20260819_22:03}
t
&=
-
s
+
\frac{
1
-
s^2
}{
2x
}.
\end{align}
Write $t(x,s)$ for $t$ as in \zcref{20260819_22:03}, and note that this is a decreasing bijection between $s \in [-1,1]$ and $t \in [-1,1]$. 

\medskip

Observe that $B_1 \cap B(\x,x+s)$ is the union of two spherical caps:
\[
\{\y =  (y_1, y_2, \ldots, y_d) \in B_1:y_1\geq t(x,s)\}
\cup
\{\y =  (y_1, y_2, \ldots, y_d) \in B(\x,x+s):y_1\leq t(x,s)\}
.
\]
\begin{figure}[htbp]
\begin{subfigure}{0.5\textwidth}
\centering
\includegraphics{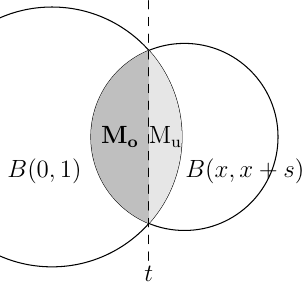}
\caption{$x+s<1$}
\end{subfigure}
\begin{subfigure}{0.49\textwidth}
\centering
\includegraphics{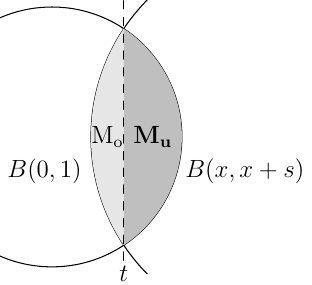}
\caption{$x+s>1$}
\end{subfigure}
\caption{The intersection of two balls is the union of two caps. It depends on the radii of the balls which of the caps is the larger one.}
\end{figure}
According to our normalization in \zcref{Sec4.1_cap}, the volume of the first cap is $\omega_d\, C(t(x,s))$, while volume of the second cap is
\(
\omega_d
(x+s)^d\,
C\bigl(
\frac{
x-t(x,s)
}{
x+s
}
\bigr)
\). We divide by the volume of $B(\x,x+s)$ and define two auxiliary maximal functions: $\Mu$ that considers the cap of the unit ball, and $\Mo$ that considers the cap of the outer ball, i.e. 
\begin{align*}
\Mu(x)
&:=
\sup_{-1\leq s\leq 1}
\frac{
C(t(x,s))
}{
(x+s)^d
} & ; & & 
\Mo(x)
&:=
\sup_{-1\leq s\leq 1}
C\left(
\frac{
x-t(x,s)
}{
x+s
}
\right).
\end{align*}
With $M(x) :=\M\chi_{B_1}(\x)$ being the maximal function of the unit ball, we plainly have the relation 
\begin{equation}\label{20260824_22:39}
M(x)  \sim \max\big\{\Mu(x), \Mo(x)\big\}.
\end{equation}
We now embark on a quest to understand who wins between $\Mu(x)$ and $\Mo(x)$ and to quantify their values. Observe that, when $x +s \leq 1$, the dominant cap is the one from the outer ball, i.e.,
\begin{align}\label{20260820_11:51}
\sup_{\substack{-1\leq s\leq 1 \\ x+s \leq 1}}
\frac{
C(t(x,s))
}{
(x+s)^d}
\leq \Mo f(x).
\end{align}

\newcommand\tlowerbound{\ensuremath{20\,d^{-\frac12}}}

\subsection{The cap of the outer ball}
\label{subsec_capbxr*}
We start by giving an asymptotic estimate for the volume of the cap of the outer ball.
\begin{lemma}\label{lemma_Mo}
With notations as above, for $x>1$ we have
\begin{equation*}
\Mo(x)
\sim
\frac{
x^{-(d+1)}
}{
\sqrt d\sqrt{1-1/x}
+
1
}
.
\end{equation*}
\end{lemma}
\begin{proof}
In order to compute $\Mo(x)$, we only need to find the minimum of
\begin{equation}\label{20260805_16:40}
\frac{
x-t(x,s)
}{
x+s
}
\end{equation}
in $s$ since $C$ is decreasing. First note that at the endpoints $s\in\{-1,1\}$ we have $(x-t(x,s))/(x+s)=1$, and so these are not points of minimum of \zcref{20260805_16:40}. That means it suffices to investigate local extrema which occur where the derivative is zero, i.e.\ at
\begin{equation}\label{20260805_16:46}
0
=
\frac{
-\partial_s t(x,s)
}{
x+s
}
-
\frac{
x-t(x,s)
}{
(x+s)^2
}
.
\end{equation}
Since $\partial_s t(x,s) = -1 -\frac sx = -\frac{x+s}x$, recalling \zcref{20260806_16:18} and \zcref{20260819_22:03} we find that equation \zcref{20260805_16:46} has the unique solution 
\begin{align}\label{20260820_09:17}
s = \sqrt{x^2-1} -x  \ \ \ ; \ \ \ t(x,s) = \frac{1}{x}.
\end{align}
From \zcref{pro_capvolume}, 
\[
\Mo(x)
=
C\left(
\frac{
x-t(x,s)
}{
x+s
}\right)
=
C\left(\sqrt{1-1/x^2}\right)
\sim
\frac{
x^{-(d+1)}
}{
\sqrt d\sqrt{1-1/x^2}
+
1
}
\sim
\frac{
x^{-(d+1)}
}{
\sqrt d\sqrt{1-1/x}
+
1
}
.
\]
\end{proof}
As a by-product of the proof above, we have the following result.
\begin{lemma}
For $x\geq \sqrt{2}$ we have $\Mo(x) \leq \Mu(x)$. 
\end{lemma}
\begin{proof}
From \zcref{20260820_09:17} in the proof above, the maximum average defining $\Mo$ occurs when $x+s = \sqrt{x^2 -1} \geq 1$. For this particular $s$, the corresponding average on the unit ball already beats $\Mo(x)$ and the result follows. 
\end{proof}

\subsection{The cap of the unit ball}
We now consider the function $\Mu$, which is more difficult to estimate. For $-1 \leq t \leq 1$ let
\[
D(t)
:=
\begin{cases}
1
&
-1\leq t\leq d^{-\frac12}
\\
\frac{
(1-t^2)^{\frac{d+1}2}
}{
t\sqrt d
}
&
d^{-\frac12}\leq t\leq 1
.
\end{cases}
\]
Recall by \zcref{pro_capvolume} that \(C(t)\sim D(t)\) (note that this holds also the range $-1 \leq t \leq 0$).
We have to maximize $C(t(x,s))/(x+s)^d$, which is comparable to $D(t(x,s))/(x+s)^d$, so it suffices to estimate the maximum of the latter, i.e.,
\begin{equation}\label{20260824_11:50}
\Mu(x) \sim \sup_{-1\leq s\leq 1} \frac{D(t(x,s))}{(x+s)^d}.
\end{equation}
In the argument below, we investigate the supremum in \zcref{20260824_11:50} in two different regimes: when $s_0 \leq s \leq 1$ and when $-1 \leq s \leq s_0$. We choose the value of $s_0$ so that 
\begin{align}\label{20260825_11:23}
t_0 = t(s_0, x) = \tlowerbound.
\end{align}
Throughout the rest of the section we may assume that $d \geq d_0$ for a fixed and large $d_0$, since \zcref{20260321_15:27} already yields \zcref{Thm4} for small dimensions. Take, for instance, $d_0 = 10^6$. 

\subsubsection{Estimates for small \texorpdfstring{$t$}t (large \texorpdfstring{$s$}s)}
We start by analyzing the range in which $t$ is close to the lower endpoint ($s$ is close to the upper endpoint). 
\begin{lemma}\label{Lemma_large_t}
With notations as above, for $s_0 \leq s \leq 1$ we have
\begin{align*}
\frac{D(t(x,s))}{(x+s)^d} \lesssim \frac{D(t(x,s_0))}{(x+s_0)^d} 
\sim
\frac{
x^{-d}
}{
(
1
-
40 /(\sqrt d\, x)
+
1/x^2
)^{\frac d2}
}
.
\end{align*}
\end{lemma}
\begin{proof}
Since $(x+s)^d$ is minimal if $t$ is maximal, and $D(t)$ remains approximately constant for $-1\leq t\leq d^{-1/2}$, it suffices to consider $t\geq d^{-\frac12}$. By \zcref{lem_capvolumedisturbance}, one observes that the values of $D(t)$ are comparable for $t$ in the range $d^{-1/2} \leq t \leq \tlowerbound$ (in fact we actually apply Lemma \ref{lem_capvolumedisturbance} five times, at each time doubling the size of the entry of $D$, and note that the exponential decay factor appearing in \zcref{20260805_14:09} is $O(1)$ at this scale $t\sim d^{-\frac12}$). Hence, it suffices to consider matters at the particular point $t_0 = t(s_0,x) = \tlowerbound$. Recall from \zcref{20260806_16:18} that
\[
(x+s_0)^2
= x^2
-2t_0x
+1 \,,
\]
and hence, 
\begin{align*}
\frac{D(t(x,s_0))}{(x+s_0)^d}
\sim
\frac1{
(
x^2
-2t_0x
+1
)^{\frac d2}
}
 =
\frac1{
(
x^2
-40x/\sqrt d
+1
)^{\frac d2}
}
=
\frac{
x^{-d}
}{
(
1
-
40 /(\sqrt d\, x)
+
1/x^2
)^{\frac d2}
}
.
\end{align*}
\end{proof}

\subsubsection{Estimates for  large \texorpdfstring{$t$}t (small \texorpdfstring{$s$}s)} 

The case $s=-1$ yields $t = 1$ and $D(1)=0$. We are then left with the task of estimating the local maxima of $s \mapsto D(t(x,s))/(x+s)^d$, when $-1 < s < s_0$, i.e. when $\tlowerbound < t(x,s) <1$. These occur when the derivative is zero. To ease the reading, we write $t$ for $t(x,s)$, and consider
\begin{align}\label{20260825_11:25}
\begin{split}
\frac\d{\d s} \frac{ D(t)}{(x+s)^d} & =  \left(
\frac{
-d
}{
x+s
}
+
\frac{
D'(t)\partial_s t
}{
D(t)
}
\right)
\frac{
D(t)
}{
(x+s)^d
}\\
& = \left(\left( \frac{-d}{x+s} \right)+ \left(
\frac{
dt^2+1
}{
(1-t^2)t
}
\right)\left(\frac{x+s}{x}\right) \right)  \frac{ D(t)}{(x+s)^d}\\
& = \left( -1 + \left(
\frac{
t^2+\frac{1}{d}
}{
(1-t^2)t
}
\right)\left(\frac{(x+s)^2}{x}\right) \right)  \frac{ d\, D(t)}{(x+s)^{d+1}}.
\end{split}
\end{align}
Hence
\begin{align*}
\frac\d{\d s} \frac{ D(t)}{(x+s)^d} = 0 \ \iff \  \frac{
(1-t^2)t
}{
t^2+1/d
}
=
\frac{
(x+s)^2
}{
x
}
=
\frac{
(x-t)^2
+
1-t^2
}{
x
}\,,
\end{align*}
where we have used \zcref{20260806_16:18} in the last identity. Defining
\[
F(x,t,1/d)
:=
\frac{
(1-t^2)t
}{
t^2+1/d
}
-
\frac{
(x-t)^2
+
1-t^2
}{
x
}\,,
\]
what we need to do is solve $F(x,t,1/d)=0$. Since it is not easy to solve $F(x,t,1/d)=0$ for $t$ explicitly (a cubic equation in $t$), our strategy here is to first solve $F(x,t,0)=0$ (i.e. replacing the parameter $1/d$ by $0$) and compute $D(t)/(x+s)^d$ for its solution. Then we will we show that this value is actually comparable to the value associated to the solution of $F(x,t,1/d)=0$.
We do so by first obtaining rough estimates for the solution of $F(x,t,1/d)=0$.
Next, in this range we prove an upper bound on $|F(x,t,1/d)-F(x,t,0)|$ and a lower bound on $|\partial_t F(x,t,0)|$.
This means that any solution of $F(x,t,1/d)=0$ has a solution of $F(x,t,0)=0$ close by, in a suitable scale for our purposes.

\begin{lemma}\label{Lem_Sol_0}
If $\tlowerbound < t <1$ and if $F(x,t,0)=0$, then $ t = 1/x$ and 
\begin{align*}
\frac{D(t)}{(x+s)^d} \sim 
\frac{
x^{-d+1}\sqrt{1-1/x}
}{
\sqrt d
}
.
\end{align*}
\end{lemma}
\begin{proof}
A standard computation yields that the only admissible solution of 
$$0 = F(x,t,0) = \frac{
(1-t^2)
}{
t
}
-
\frac{
(x-t)^2
+
1-t^2
}{
x
}$$
is $t = 1/x$ (note that this happens to be the same solution as for $\Mo$). In our range, that means that $x < \sqrt{d}/{20}$ and thus
\begin{align*}
\frac{D(t)}{(x+s)^d}  & = \frac{D(t)}{((x-t)^2+1-t^2)^{d/2}}
=
\frac{D(1/x)}{(x^2-1)^{\frac d2}}
=
\frac x{\sqrt d}
\frac{
(1-1/x^2)^{\frac{d+1}2}
}{
(x^2-1)^{\frac d2}
}\\
& =
\frac{x^{-d} \sqrt{x^2-1}}{\sqrt{d}}
=
\frac{
x^{-d+1}\sqrt{1-1/x^2}
}{
\sqrt d
}
\sim
\frac{
x^{-d+1} \sqrt{1-1/x}
}{
\sqrt d
}
.
\end{align*}
\end{proof}
We now make a simple, yet very important, assumption. Recall that if $x + s \leq 1$, we have already observed in \zcref{20260820_11:51} that $D(t)/(x+s)^d \lesssim \Mo(x)$. Therefore, going forward, we may restrict our attention to the investigation of the local extrema under the following conditions:
\begin{align}\label{20260807_13:49}
\tlowerbound
<
t<1 \ \ \ \ \ ; \ \ \ \ \ 
x
>
1 \ \ \ \ \ ; \ \ \ \ \ 
(x-t)^2+1-t^2
=
(x+s)^2
>
1
.
\end{align}

\begin{lemma}[A priori bounds on $t$] \label{Lemma_a_priori_bounds}
Under conditions \zcref{20260807_13:49}, if $F(x,t,1/d)=0$ then 
\begin{equation}\label{20260824_13:12}
\frac1x
\geq
t
\geq
\frac1{12\, x}\,.
\end{equation}
\end{lemma}
\begin{proof} The first inequality in \zcref{20260824_13:12} holds in general. In fact, if $F(x,t,1/d)=0$ then 
\begin{equation}\label{20260825_11:56}
\frac{1 - t^2}{t} \geq \frac{
(1-t^2)t
}{
t^2+1/d
} = \frac{
x^2 - 2tx +1
}{
x}\,,
\end{equation}
which yields
\begin{align}\label{20260825_11:57}
x - xt^2 \geq x^2t - 2t^2x + t \ \ \Rightarrow \ \ (x-t)(1 - xt)\geq0\ \ \Rightarrow \ \ \frac{1}{x} \geq  t.
\end{align}
For the second inequality in \zcref{20260824_13:12} we shall use conditions \zcref{20260807_13:49}. Since $x+s \geq 1$, and $s \geq -1$ we have $(x+s)^2/x\geq 1/2$.
Thus, $F(x,t,1/d)=0$ and $t\geq0$ require
\begin{align}\label{20260807_16:17}
\frac{
1-t^2
}t
&\geq
\frac{
(1-t^2)t
}{
t^2+1/d} = \frac{
(x+s)^2
}{
x
}
\geq
\frac12  \ \ \Rightarrow  \ \ t \leq \frac{\sqrt{17} -1}{4} < \frac{4}{5}.
\end{align}
Under \zcref{20260807_13:49}, observe also that 
\begin{align}\label{20260807_16:18}
\frac{
(x+s)^2
}x
\leq
4x
.
\end{align}
From \zcref{20260807_16:17} we get $ \frac{9}{25} \leq1-t^2$. Using this, combined with \zcref{20260807_16:18} and the equation  $F(x,t,1/d)=0$ as in \zcref{20260807_16:17}, we get 
\begin{align*}
\frac t{t^2+1/d} = \frac{(x+s)^2}{x (1 - t^2)} 
\leq \tfrac{100}{9} x
 \ \ \Rightarrow \ \ 
t
+\frac1{dt}
\geq
\frac{1}{\frac{100}{9} x}
.
\end{align*}
Since $t\geq\tlowerbound$ we have
\(
t+\frac1{dt}
\leq
\big( 1+ \frac{1}{400}\big)t
\). Hence,
\[
t
\geq \frac{t+\frac1{dt}}{1+ \frac{1}{400}}  \geq \frac{1}{\left(1+ \frac{1}{400}\right)}\frac{1}{\tfrac{100}{9}x} \geq 
\frac1{12 x}
.
\]
\end{proof}

\begin{lemma}[Comparison lemma]\label{Lemma_comparison}
Under conditions \zcref{20260807_13:49}, if $F(x,t,1/d)=0$ there exists $\tau$ with $t \leq \tau \leq  \frac{11}{10}t$ and $(1-\tau^2)\geq(1-t^2)/2$ verifying $F(x,\tau,0)=0$ and such that 
 \begin{align*}
 0 \leq \tau - t \leq16 \,\varepsilon(t)\,,
 \end{align*}
 with $\varepsilon(t)$ as in \zcref{20260807_16:53}. Furthermore,
 \begin{align}\label{20260811_17:33}
\frac{D(t)}{(x+s)^d} \sim 
\frac{
x^{-d+1}\sqrt{1-1/x}
}{
\sqrt d
}
.
 \end{align} 
\end{lemma}
\begin{proof}
The point of this lemma is to control the dependence of $F(x,t,1/d)$ on the parameter $1/d$. Start by recalling that 
\begin{align*}
F(x,t,0)
&=
\frac1t
-t
-
\frac{
(x-t)^2
+
1-t^2
}x = \frac{1}{t} + t - \frac{x^2 +1}{x}.
\end{align*}
We compute derivatives in $t$,
\begin{align}\label{202600807_17:10}
F'(x,t,0)
&=
-\frac1{t^2}
+1
=
\frac{t^2-1}{t^2}
.
\end{align}
Note also that 
\begin{align}\label{202600807_17:11}
F(x,t,0)
-
F(x,t,1/d)
=
\frac{
1-t^2
}{
t(dt^2+1)
}.
\end{align}
Let $\tau$ with $t \leq \tau \leq  \frac{11}{10}t$ and $1-\tau^2\geq(1-t^2)/2$.
Then, from \zcref{202600807_17:10}, \zcref{202600807_17:11} and the fact that $t \geq \tlowerbound$, 
\begin{align} \label{20260811_14:19}
\begin{split}
&\frac{
F(x,t,1/d)
 -
F(x,t,0)
}{
F'(x,\tau,0)
}
=
\frac{(1-t^2)}{(1-\tau^2)}
\frac{\tau^2}{t^2}
\frac{t}{
(dt^2+1)
}
\\
& \ \ \ \ \ \ \ \ \ \ \ \ \  \ \ \leq 2 \left(\frac{11}{10}\right)^2 \frac{t}{dt^2 +1} \leq
\frac{
3t
}{
dt^2+1
}
=
\frac{
3
}{
dt+1/t
}
\leq
\frac{
3
}{
d\cdot \tlowerbound
}
\leq
\frac t{10}
.
\end{split}
\end{align}
Now, assume $F(x,t,1/d)=0$.
Then, from \zcref{20260807_16:17} we already know that $t\leq\frac45$.
In particular, if $0\leq \tau-t\leq t/10\leq1/10$, then
\begin{align*}
2(1-\tau)
-
(1-t)
=
1-2(\tau - t) - t
\geq
1-1/5-t
\geq0
\end{align*}
and hence
\begin{align*}
\frac{1-\tau^2}{1-t^2}
&=
\frac{(1-\tau)(1+\tau)}{(1-t)(1+t)}
\geq
\frac{1-\tau}{1-t}
\geq
\frac12
.
\end{align*}
From \zcref{202600807_17:11} note that $F(x,t,0) >0$, and from \zcref{202600807_17:10} we have $F'(x,\tau,0) <0$.  From the fundamental theorem of calculus and \zcref{20260811_14:19}, 
\begin{align*}
F(x,t + \tfrac{t}{10},0) & = F(x,t,0) + \int_{t}^{t + \frac{t}{10}} F'(x,\tau,0) \, \d\tau\\
& \leq  F(x,t,0) \left( 1-   \frac{10}{t}   \int_{t}^{t + \frac{t}{10}}   \d\tau\right) = 0.
\end{align*} 
By the mean value theorem, there exists $t< \tau \leq\frac{11}{10}t$ with $F(x,\tau,0)=0$, which, recall, means $\tau=1/x$. Another application of the mean value theorem shows that there exists $\tilde{\tau}$ with $t \leq  \tilde{\tau} \leq \tau$ such that 
\begin{align}\label{20260811_14:51}
0 \geq F'(x, \tilde{\tau},0)  \geq \frac{F(x, \tau ,0)- F(x,t,0)}{\tau - t}.
\end{align}
Since $t \leq \frac45$ we have
\begin{align}\label{20260811_16:36}
1 - \tilde{\tau}^2 \geq 1 - \tau^2 \geq 1 - \left(\frac{11t}{10}\right)^2 \geq 1 - \left(\frac{22}{25}\right)^2 > \frac{1}{5}.
\end{align}
Using \zcref{20260811_14:51} and \zcref{20260811_16:36}, and recalling that $t\geq\tlowerbound$, we get 
\begin{align} \label{20260811_16:53}
\begin{split}
0 \leq
\tau-t
& \leq \frac{F(x, \tau ,0)- F(x,t,0)}{F'(x, \tilde{\tau},0)} = 
\frac{F(x,t,1/d)
-
F(x,t,0)
}{
F'(x,\tilde{\tau},0)
}
\\
& = \frac{(1-t^2)}{(1-\tilde{\tau}^2)}
\frac{\tilde{\tau}^2}{t^2}
\frac{t}{
(dt^2+1)
}\\
& \leq
5 \left(\frac{11}{10}\right)^2 \frac{
(1-t^2)
}{
dt+1/t
}
\leq
\frac{
7(1-t^2)
}{
dt
}
\leq
\frac{
8(1-t^2)
}{
dt+\sqrt d
}
=
16\varepsilon(t)
.
\end{split}
\end{align}
An application of \zcref{lem_capvolumedisturbance} yields that $D(t) \sim D(\tau)$. It remains to compare $(x + s(x,t))^2$ and $(x + s(x,\tau))^2$. By \zcref{20260806_16:18} these are equal, respectively, to $(x-t)^2+1-t^2$ and $(x-\tau)^2+1-\tau^2$. Since $t \leq \tau$ we have
\[
(x-\tau)^2+1-\tau^2
\leq
(x-t)^2+1-t^2
.
\]
For the reverse inequality, we start by recalling \zcref{20260811_16:53} and the fact that $t \geq 1/(12 x)$ by \zcref{Lemma_a_priori_bounds} to get
\[
0 \leq \tau-t
\leq
\frac{
7(1-t^2)
}{
dt
}
\leq
\frac{7}{dt} \leq 
\frac{84x}d
.
\]
That, in turn, yields (recall $\tau \leq \frac{11t}{10} \leq \frac{22}{25}$)
\begin{align*}
\big((x-t)^2+1-t^2\big)
-
\big((x-\tau)^2+1-\tau^2\big) = 2x(\tau - t) 
&\lesssim \frac{x^2}{d}
\lesssim
\frac{
(x-\tau)^2+1-\tau^2
}d.
\end{align*}
Hence,
\begin{align*}
1 \leq \left(\frac{
(x-t)^2+1-t^2
}{
(x-\tau)^2+1-\tau^2}\right)^{\frac d2}
 = \left(1 + \frac{\big((x-t)^2+1-t^2\big)
-
\big((x-\tau)^2+1-\tau^2\big)}{(x-\tau)^2+1-\tau^2}\right)^{\frac d2}&\lesssim
1
.
\end{align*}
We then arrive at the conclusion that 
\begin{align*}
\frac{D(t)}{(x+s(x,t))^d} \sim \frac{D(\tau)}{(x+s(x,\tau))^d} \,,
\end{align*}
and \zcref{20260811_17:33} plainly follows from \zcref{Lem_Sol_0}.

\end{proof}

\subsection{Conclusion: proof of \texorpdfstring{\zcref{Thm4}}{Theorem 1.4}}
Recall that we may assume here that $d \geq 10^6$ since \zcref{20260321_15:27} yields \zcref{Thm4} for small dimensions. From \zcref{20260824_22:39}, \zcref{lemma_Mo}, \zcref{Lemma_large_t}, and \zcref{Lemma_comparison}, we have
\begin{align*}
\M f(x) \sim x^{-d} \cdot \max\left\{ \frac{
x^{-1}
}{
\sqrt d\sqrt{1-1/x}
+
1
} \ ,\   \frac{
1
}{
(
1
-
40 /(\sqrt d\, x)
+
1/x^2
)^{\frac d2}
} \ ,\  \frac{
x\sqrt{1-1/x}
}{
\sqrt d
}
\right\}\,,
\end{align*}
where the last term is only admissible if $s \mapsto D(t(x,s))/(x+s)^d$ has  a local maximum at a point $s$ verifying conditions \zcref{20260807_13:49}. We now carefully analyze which term wins (here, winning means that the other two are bounded by a universal constant times it).

{\it Case 1.} The range $1 < x \leq 2$. Here the first term (the one corresponding to $\Mo f$) in the maximum wins.

{\it Case 2.} The range $\sqrt{d} \leq x$. Here we the conclusion of \zcref{Lemma_a_priori_bounds} is not verified, and hence the last term is not admissible. We then note that the second term wins since it is $\sim 1$ in this range. In fact,
\begin{align} \label{20260824_23:22}
1 \lesssim \big(1 - \tfrac{40}{ d} \big)^{\frac d2} \leq \big(1 - \tfrac{40}{\sqrt d\, x} \big)^{\frac d2}\leq \big(1 - \tfrac{40}{\sqrt d\, x} + \tfrac{1}{x^2}\big)^{\frac d2} \leq 1.
\end{align}

{\it Case 3.} The range $2 \leq x \leq \sqrt{d}$. We further subdivide this range in two. 

If $ \sqrt{d} /80 \leq x \leq \sqrt{d} $, a similar computation as in \zcref{20260824_23:22} shows that both the second and third terms are $\sim 1$ (and hence win against the first). 

If $2\leq x \leq  \sqrt{d} /80$ we claim that the third term is admissible and wins. That it wins against the first one is clear. To see that it wins against the second one, note that for a universal $C>0$ we have (in the second inequality below we use $\log(1 +y) \geq y/2$ for $0 \leq y \leq 1$)
\begin{align*}
& \frac{1}{2} \log \left(\frac{d}{x^2} \right) \leq C + \frac{d}{8x^2} \leq C + \frac{d}{2} \log\left(1+\frac{1}{2 x^2}\right)  \Rightarrow  \frac{\sqrt{d}}{x} \lesssim  \left(1+\frac{1}{ 2 x^2}\right)^{\frac d2}\\
& \Rightarrow \frac{x}{\sqrt{d}} \gtrsim \frac{1}{\left(1+\frac{1}{2 x^2}\right)^{\frac d2}} \geq \frac{1}{\left(1-\frac{40}{\sqrt d\, x} + \frac{1}{x^2}\right)^{\frac d2}}.
\end{align*}
It remains to verify that $s \mapsto D(t(x,s))/(x+s)^d$ has  a local maximum at a point $s$ verifying conditions \zcref{20260807_13:49}. Recall the definition of $t_0$ and $s_0$ in \zcref{20260825_11:23}. From \zcref{20260825_11:25} we have
\begin{align}\label{20260825_11:31}
\frac\d{\d s} \frac{ D(t)}{(x+s)^d} < 0 \iff 
\frac{(1-t_0^2)t_0}{t_0^2+1/d}>\frac{(x+s_0)^2}{x}=\frac{(x-t_0)^2+1-t_0^2}{x} = \frac{x^2 - 2 x t_0 +1}{x}.
\end{align}
Since $t_0= \tlowerbound$, the validity of \zcref{20260825_11:31} follows from the fact that 
\begin{align*}
\frac{
(1-t_0^2)t_0
}{
t_0^2+1/d
}   = \frac{\left(1 - \frac{400}{d}\right) \frac{20}{\sqrt{d}}}{\frac{401}{d}}   = \frac{20\sqrt{d} -  \frac{8000}{\sqrt{d} }}{401}  > \frac{\sqrt{d}}{40} \geq 2x \geq \frac{x^2+1}{x} \geq  \frac{x^2 - 2 x t_0 +1}{x}
\end{align*}
in our range. This shows indeed that there exists a maximum of $s \mapsto D(t(x,s))/(x+s)^d$ in the interval $-1 < s < s_0$ (or, equivalently, $t_0 < t < 1$), and that this maximum beats the value of the function at the particular point $t_0= \tlowerbound$. Recall that such a point of maximum $t$ must solve 
$$\frac{(1-t^2)t}{t^2+1/d} = \frac{x^2 - 2tx +1}{x}.$$
In this case, we have seen from \zcref{20260825_11:56} and \zcref{20260825_11:57} that $t \leq \frac{1}{x}$. Since $x \geq 2$ we have  $\frac{1}{x} < \frac{x}{2}$. Hence
\begin{align*}
1 < x^2 - 2tx +1  = (x-t)^2+1-t^2 = (x-s)^2\,,
\end{align*}
and this verifies indeed \zcref{20260807_13:49} and completes the proof.

\section{Proof of \texorpdfstring{\zcref{Cor5}}{Corollary 1.5}}\label{Sec5_Cor1.5}
Let $g(\x) = \rad{g}(|\x|)$ be the radial function appearing on the right-hand side of \zcref{20260912_16:07}. Since this is defined in \zcref{20260912_16:07} just for $|\x|>1$, let us set here $g(\x) = \rad g(|\x|) = 1$ for $0 \leq |\x| \leq 1$. Since $\M\chi_{B_1}(\x) \sim g(\x)$, by \zcref{radiallydecreasingcomparison} derived from the coarea formula we have 
$$\var(\M\chi_{B_1}) \sim \var(g).$$
To compute $\var(g)$ we appeal to \zcref{varfradial}. Note that $\rad g$ is continuous and 
\begin{align}\label{20260912_17:12}
-x^{d-1} \, \rad{g}\, '(x)  
&=
\begin{cases}
0\, & {\rm if} \  \ \ 0 < x < 1\,;\\ \\
d \,x^{-2}\, & {\rm if} \  \ \ 1 < x \leq 1 + \frac{1}{d}\,;\\ \\
 \frac{d \, x^{-2} (x-1)^{-1/2}}{\sqrt{d}} +  \frac{x^{-1} (x-1)^{-3/2}}{2\sqrt{d}} \, , & {\rm if} \  \ 1 + \frac{1}{d} < x < 2\,;\\
\\
 \frac{(d-1) x^{-1}}{2 \sqrt{d}}\, ,  & {\rm if} \ \ 2 < x < 2\sqrt{d}\,;\\
\\
d \,x^{-2} \, , & {\rm if} \  \ 2\sqrt{d} < x.
\end{cases}
\end{align}
Then, from \zcref{varfradial,20260912_17:12}, (recall that $\var(\chi_{B_1}) = \sigma_{d-1}$)
\begin{align*}
\frac{\var(\M\chi_{B_1})}{\var(\chi_{B_1})}  \sim \frac{\var(g)}{\sigma_{d-1}} = 
-\int_0^\infty
x^{d-1} \, \rad g\,'(x)  \d x \sim  \sqrt{d}\,\log d.
\end{align*}
The main term comes from the interval $2 < x < 2\sqrt{d}$, while the other intervals contribute $O(\sqrt{d})$. This completes the proof of \zcref{Cor5}.

\begin{remark}\label{Remark_coarea_HL}
Observe that if we just use the bounds in \zcref{20260912_17:58} we only get, via \zcref{radiallydecreasingcomparison} derived from the coarea formula and \zcref{varfradial}, 
\begin{align*}
2^d-1 =  d \int_{2}^{\infty} \frac{x^{d-1}}{(x-1)^{d+1}}\,\d x \geq   \frac{\var(\M\chi_{B_1})}{\sigma_{d-1}} \geq  d \int_{1}^{\infty} \frac{x^{d-1}}{(x+1)^{d+1}}\,\d x = 1 - 2^{-d}.
\end{align*}
\end{remark}

\section{Proof of \texorpdfstring{\zcref{Cor6}}{Corollary 1.7}}

The result is true when $d=1$, so we may assume that $d \geq 2$. Let $f \in \BV(\R^d)$ be non-negative and radially decreasing. We may assume that $C = \lim_{|\x| \to \infty} f(\x) = 0$, otherwise we work with the function $f_1:= f - C$. Then, by the layer cake representation,
\[
f(\x) = \int_{0}^{\infty} \chi_{\{f >t\}} (\x) \d t
,
\]
and recall the coarea formula (\zcref{coarea}).
Note that the set $\{\x\in\mathbb{R}^d:f(\x)>t\}$ is a ball.
Let us write $g_t(\x):= \M\chi_{\{f >t\}} (\x)$ and
\[
G(\x) :=  \int_{0}^{\infty} g_t(\x) \d t\,,
\]
and write $G(\x) = \rad G(|\x|)$.
Recall from \zcref{20260321_15:27} that 
\[
g_t(\x) \sim_d
\begin{cases}
1,&\text{if }f(\vt x)> t
,\\
\frac{|\{f >t\}|}{\omega_d |\x|^d},&\text{if }f(\vt x)\leq t
.
\end{cases}
\]
Thus, using the coarea formula we get 
\begin{align*}
G(\x) & \lesssim_d \int_0^{f(\x)} 1 \d t + \frac{1}{\omega_d |\x|^d}\int_{f(\x)}^{\infty} |\{f >t\}| \d t = f(\x) + \frac{1}{\omega_d |\x|^d}\int_{f(\x)}^{\infty} |\{f >t\}|^{\frac{d-1}{d}}\,  |\{f >t\}|^{\frac{1}{d}} \d t\\
&  \leq f(\x) + \frac{\omega_d^{1\d} |\x|}{\omega_d |\x|^d}\int_{f(\x)}^{\infty} |\{f >t\}|^{\frac{d-1}{d}}\d t  \lesssim_d f(\x) + \frac{1}{|\x|^{d-1}} \var (f). 
\end{align*}
The last term is finite for any $\vt x\neq0$, integrates on any ball which means $G\in L^1_\loc(\mathbb{R}^d)$, and moreover tends to $0$ as $|\x| \to \infty$.

Then \zcref{20260903_10:40} and the sublinearity of the Hardy-Littlewood maximal operator imply
\[
\M_{\varphi}f(\x) \leq \M f(\x)
\leq
G(\vt x)
.
\]
Since $\M_\varphi f$ and $G$ are radial decreasing functions with $\lim_{|\x| \to \infty} G(\x) = 0$, by \zcref{radiallydecreasingcomparison} it suffices to prove that 
\[
\var (G)
\lesssim
\sqrt d\log d\,
\var (f)
.
\]
By \zcref{gradfdisintegrationformulaweak}, if $\varphi \in C^{1}_\cpt((0,\infty))$ with $|\varphi(r)|\leq r^{d-1}$, we have (the use of Fubini's theorem below is justified since $G\in L^1_\loc(\mathbb{R}^d)$)
\begin{align*}
\sigma_{d-1}  \int_0^\infty \rad G(r) \,\varphi'(r)\,\d r & = \sigma_{d-1} \int_0^{\infty}\left( \int_{0}^{\infty} \rad g_t(r) \d t\right) \varphi'(r) \d r =  \sigma_{d-1}  \int_0^\infty \int_0^\infty \rad g_t(r) \, \varphi'(r) \d r \d t \\
& \leq \int_0^\infty
\var(g_t)
\intd t.
\end{align*}
Taking the supremum over $\varphi$, applying \zcref{approx_by_annuli,gradfdisintegrationformulaweak} yields
\begin{equation}
\label{derivative_inside_integral}
\var (G) \leq \int_0^\infty
\var(g_t)
\intd t.
\end{equation}
Applying \zcref{Cor5} and the coarea formula we obtain
\[
\int_0^\infty
\var(g_t)
\intd t
\lesssim
\sqrt d\log d
\int_0^\infty
\var(\chi_{\{f>t\}})
\intd t
=
\sqrt d\log d\,
\var (f)
<\infty
,
\]
which concludes the proof.

\begin{remark}
With some extra work, it is actually possible to show equality in \zcref{derivative_inside_integral}, but this is not strictly needed for our purposes. 
\end{remark}

\appendix

\section{Radon measures and functions of bounded variation}
\label{appendix_radonbv}

Given a Radon measure $\mu$ on $\mathbb{R}^d$ we denote by $|\mu|$ its total variation measure and by $\|\mu\|$ its total variation.

\subsection{Functions of bounded variation}

For the following \zcref*[noref,nocap]{definition_bv} we follow \cite[Section~3.1]{AmbrosioFuscoPallara}.

\begin{definition}
\label{definition_bv}
Let $\Omega\subset\mathbb{R}^d$ be open.
We denote by $\BV(\Omega)$ the set of all $f\in L^1_\loc(\Omega)$ whose gradient can be represented as a Radon measure, in the sense that there exists an $\R^d$-valued Radon measure $\grad f$ on $\Omega$ such that for all $\vtf \varphi\in C^1_\cpt(\Omega,\mathbb{R}^d)$ we have
\[
\int_{\Omega}
f \Div \vtf \varphi 
=
\int_{\Omega}
 \vtf \varphi \cdot \d \grad f
.
\]
Moreover, for $f\in L^1_\loc(\Omega)$ denote
\begin{equation}\label{20261008_21:50}
\var_\Omega(f)
=
\sup\Bigl\{
\int_{\Omega}
f \Div \vtf \varphi 
:
\vtf \varphi \in C^1_\cpt(\Omega,\mathbb{R}^d)
,\ 
\|\vtf \varphi\|_{L^\infty(\Omega)}\leq 1
\Bigr\}
.
\end{equation}
Then $f\in\BV(\Omega)$ if and only if $\var_\Omega(f)<\infty$, in which case $|\grad f|(\Omega)=\|\grad f\|=\var_\Omega(f)$.

We denote by $\BV_\loc(\Omega)$ the set of all $f\in L^1_\loc(\Omega)$ such that, for all open $\Omega'  \subset  \subset \Omega$, we have $f\in\BV(\Omega')$.
\end{definition}

For the following coarea formula we refer the reader to \cite[Section~3.3]{AmbrosioFuscoPallara}.
For a measurable set $E\subset\Omega$ we denote its perimeter $P(E,\Omega)=\var_\Omega(\chi_E)$.
If $E$ has a $C^1$-boundary then $P(E,\Omega)=\sm{\Omega\cap\partial E}$.

\begin{lemma}[Coarea formula]
\label{coarea}
Let $f\in L^1_\loc(\Omega)$.
Then
\[
\var_\Omega(f)
=
\int_{-\infty}^\infty
P(\{\x\in\mathbb{R}^d:f(\x)>t\},\Omega)
\d t
.
\]
\end{lemma}

It has the following consequence:

\begin{lemma}
\label{radiallydecreasingcomparison}
Let $f,g:\mathbb{R}^d\rightarrow[0,\infty)$ be radially decreasing such that $g\leq f$ and $\lim_{|\x|\rightarrow\infty}f(\x)=0$.
Then
\[
\var_{\mathbb{R}^d}(g)
\leq
\var_{\mathbb{R}^d}(f)
.
\]
\end{lemma}

\begin{proof}
For any $t>0$ both $f^{-1}((t,\infty))$ and $g^{-1}((t,\infty))$ are (possibly empty) balls with \(g^{-1}((t,\infty))\subset f^{-1}((t,\infty))\).
Thus \(P(\{\x\in\mathbb{R}^d:g(\x)>t\},\mathbb{R}^d)\leq P(\{\x\in\mathbb{R}^d:f(\x)>t\},\mathbb{R}^d)\) and the claim follows from \zcref{coarea}.
\end{proof}

Denote an annulus by $A(a,b)=\{\vt x:a<|\vt x|< b\}$.
We also need the following approximation of $\R^d$ by annuli.

\begin{lemma}
\label{approx_by_annuli}
Let $f:\R^d\rightarrow\R$, $C>0$ and $a_n\rightarrow0$, $b_n\rightarrow\infty$ such that $f\in L^1_\loc(A(a_n,b_n))$ with $\var_{A(a_n,b_n)}(f)\leq C$.
Then $f\in L^1_\loc(\mathbb{R}^d)$ with $\var_{\R^d}(f)\leq C$.
\end{lemma}

This should be standard, but for the convenience of the reader we give a proof since we have not been able to find one in the literature.

\begin{proof}
\emph{Step 1:}
Let $R>0$.
We show that $f\in L^{\frac d{d-1}}(B(\vt 0,R))$.

Denote the integral average by
\(
f_{A(a,b)}
=
|A(a,b)|^{-1}\int_{A(a,b)}f
.
\)
Recall the Poincaré inequality on John domains, for example \cite[Theorem~107]{lecturenoteshajlasz}.
It applies in particular to sufficiently thick annuli, i.e.\ for $0\leq a< b<\infty$ with $b-a\gtrsim b$, where it asserts
\begin{equation}
\label{poincare_annulus}
\|f-f_{A(a,b)}\|_{L^{\frac d{d-1}}(A(a,b))}
\lesssim_d
\var_{A(a,b)}(f)
.
\end{equation}
As a consequence,
\begin{align*}
|
f_{A(r,2r)}
-
f_{A(r/2,2r)}
|
&\sim_d
r^{-d}
\left|
\int_{A(r,2r)}
f
-
f_{A(r/2,2r)}
\right|
\leq
r^{-d}
\int_{A(r/2,2r)}
|
f
-
f_{A(r/2,2r)}
|
\\
&\lesssim_d
r^{1-d}
\|f-f_{A(r/2,2r)}\|_{L^{\frac d{d-1}}(A(r/2,2r))}
\lesssim_d
r^{1-d}
\var_{A(r/2,2r)}(f)
,
\end{align*}
and doing a similar calculation for
\(
|
f_{A(r/2,r)}
-
f_{A(r/2,2r)}
|
\)
we can conclude
\begin{equation}
\label{average_difference}
|
f_{A(r/2,r)}
-
f_{A(r,2r)}
|
\lesssim_d
r^{1-d}
\var_{A(r/2,2r)}(f)
.
\end{equation}
Next, decompose
\begin{align*}
&
f-f_{A(R/2,R)}
=
\sum_{n=0}^\infty
(f-f_{A(R/2,R)})\chi_{A(2^{-n-1}R,2^{-n}R)}
\\
&\qquad=
\sum_{n=0}^\infty
\left(
f-f_{A(2^{-n-1}R,2^{-n}R)}
+
\sum_{k=1}^n
f_{A(2^{-k-1}R,2^{-k}R)}
-
f_{A(2^{-k}R,2^{-k+1}R)}
\right)
\chi_{A(2^{-n-1}R,2^{-n}R)}
.
\end{align*}
We estimate the first part using \zcref{poincare_annulus},
\begin{align*}
\sum_{n=0}^\infty
\left\|
(
f-f_{A(2^{-n-1}R,2^{-n}R)}
)
\chi_{A(2^{-n-1}R,2^{-n}R)}
\right\|_{L^{\frac d{d-1}}(\mathbb{R}^d)}
\lesssim_d
\sum_{n=0}^\infty
\var_{A(2^{-n-1}R,2^{-n}R)}(f)
\leq C
\end{align*}
and the second part using \zcref{average_difference} and Minkowski's inequality,
\begin{align*}
&
\left\|
\sum_{n=1}^\infty
\left(
\sum_{k=1}^n
|
f_{A(2^{-k-1}R,2^{-k}R)}
-
f_{A(2^{-k}R,2^{-k+1}R)}
|
\right)
\chi_{A(2^{-n-1}R,2^{-n}R)}
\right\|_{L^{\frac d{d-1}}(\mathbb{R}^d)}
\\
&\quad\sim_d
\left(
\sum_{n=1}^\infty
(2^{-n}R)^d
\left(
\sum_{k=1}^n
|
f_{A(2^{-k-1}R,2^{-k}R)}
-
f_{A(2^{-k}R,2^{-k+1}R)}
|
\right)^{\frac d{d-1}}
\right)^{\frac{d-1}d}
\\
&\quad\lesssim_d
\left(
\sum_{n=1}^\infty
\left(
\sum_{k=1}^n
2^{(n-k)(1-d)}
\var_{A(2^{-k-1}R,2^{-k+1}R)}(f)
\right)^{\frac d{d-1}}
\right)^{\frac{d-1}d}
\\
&\quad\leq
\sum_{k=1}^\infty
\left(
\sum_{n=k}^\infty
\left(
2^{(n-k)(1-d)}
\var_{A(2^{-k-1}R,2^{-k+1}R)}(f)
\right)^{\frac d{d-1}}
\right)^{\frac{d-1}d}
\\
&\quad\sim_d
\sum_{k=1}^\infty
\var_{A(2^{-k-1}R,2^{-k+1}R)}(f)
\leq 2C
.
\end{align*}
Combining the last three calculations we can conclude
\[
\|f-f_{A(R/2,R)}\|_{L^{\frac d{d-1}}(B(\vt 0,R))}
\lesssim_d
C
,
\]
in particular $f\in L^1_\loc(\mathbb{R}^d)$.

\emph{Step 2:}
Let $R>0$ and $\vtf \varphi\in C^1_\cpt(\mathbb{R}^d,\mathbb{R}^d)$ with $\|\vtf \varphi\|_{L^\infty(\R^d)}\leq1$.
Then exists $R>0$ with $\vtf\varphi\in C^1_\cpt(B(\vt 0,R),\R^d)$.
Let $\eta_n$ be a smooth cutoff with $\eta_n=0$ on $B(\vt 0,2a_n)$, $\eta_n=1$ on $A(3a_n,R)$, $|\grad \eta_n|\lesssim\chi_{A(2a_n,3a_n)}/a_n$ and $\|\eta_n\|_{L^\infty(\R^d)}=1$.
Then since $f\in L^1_\loc(\mathbb{R}^d)$ and $a_n\rightarrow0$ we have
\[
\int f\eta_n\Div(\vtf \varphi)
\rightarrow
\int f\Div(\vtf \varphi)
.
\]
Moreover, since $|\grad \eta_n|\lesssim\chi_{A(2a_n,3a_n)}/a_n$ and $f\in L^{\frac d{d-1}}(B(\vt 0,R))$ we get
\[
\|
f(\grad\eta_n\cdot\vtf \varphi)
\|_{L^1(\R^d)}
\leq
\|\grad\eta_n\|_{L^d(A(2a_n,3a_n))}
\|f\|_{L^{\frac d{d-1}}(A(2a_n,3a_n))}
\|\varphi\|_{L^\infty(\R^d)}
\lesssim
\|f\|_{L^{\frac d{d-1}}(B(\vt 0,3a_n))}
\rightarrow
0
\]
as $n\rightarrow\infty$.
Since $\eta_n\vtf\varphi\in C^1_\cpt(B(\vt 0,R))$ with $\|\eta_n\vtf \varphi\|_{L^\infty(\R^r)}\leq1$, we can conclude
\[
C
\geq
\int f\Div(\eta_n\vtf \varphi)
=
\int f\eta_n\Div(\vtf \varphi)
+
\int f(\grad\eta_n\cdot\vtf \varphi)
\rightarrow
\int f\Div(\vtf \varphi)
\]
as $n\rightarrow\infty$, finishing the proof.
\end{proof}

\subsection{Variation of radial functions}
\label{radialvariationformula}

For $f:\mathbb{R}^d\rightarrow\mathbb{R}$ radial denote by $\rad f:[0,\infty)\rightarrow\mathbb{R}$ the function with $f(\vt x)=\rad f(|\vt x|)$.
The main purpose of this \zcref*[noref,nocap]{radialvariationformula} is the disintegration formula \zcref{varfradial} for the derivative of $f$.

\begin{lemma}
\label{gradfdisintegrationformulaweak}
Let $f\in L^1_\loc(\mathbb{R}^d)$ be radial and $0\leq a< b\leq\infty$.
Then
\begin{equation*}
\var_{\{\vt x:a<|\vt x|< b\}}(f)
=
\sup\left\{
\sigma_{d-1}
\int
\varphi'
\rad f
:
\varphi\in C^1_\cpt((a,b)),\ |\varphi(r)|\leq r^{d-1}
\right\}
,
\end{equation*}
where both terms may be infinite.
\end{lemma}

\begin{proof}

\newcommand\ssym\mathring

Let $\vtf \varphi\in C^1_\cpt(\{\vt x\in\mathbb{R}^d:a<|\vt x|<b\},\mathbb{R}^d)$. Denote by $\ssym{\vtf \varphi}\in C^1_\cpt(\{\vt x\in\mathbb{R}^d:a<|\vt x|<b\},\mathbb{R}^d)$ the spherical symmetrization of $\vtf \varphi$, i.e.\
\[
\ssym{\vtf \varphi}(\vt x)
=
\int_{\SO{d}}
R^{-1}\vtf \varphi(R\x)
\intd R
,
\]
where we integrate against the normalized Haar measure on $\SO{d}$, the compact topological group of real orthogonal $d \times d$ matrices with determinant $1$.
Then $\ssym{\vtf \varphi}$ is spherically symmetric in the sense that for all $\x$ and rotations $R\in\SO{d}$ we have $\ssym{\vtf \varphi}(R\x)=R\ssym{\vtf \varphi}(\x)$.
This implies, that for every $\x$, all rotations $R$ that leave $\x$ unchanged also leave $\ssym{\vtf \varphi}(\x)$ unchanged.
That means $\ssym{\vtf \varphi}(\x)$ is parallel to $\x$, and so there exists $\varphi \in C^1_\cpt((a,b))$ with
\(
\ssym{\vtf \varphi}(\x)= \varphi(|\x|)\x/|\x|
.
\)
Let $\vtf \varphi_{\rho}$ be the radial part of $\vtf \varphi$, i.e.,\
\(
\vtf \varphi_{\rho}(\x)
=
(\vtf \varphi(\x)\cdot \x)
\frac{\x}{|\x|^2}
\).
We then have $\ssym{\vtf \varphi_{\rho}}=\ssym{\vtf \varphi}=(\ssym{\vtf \varphi})_{\rho}$,  
\[
\Div(\vtf \varphi_{\rho})(x)
=
\frac1{|\x|^d}
\frac\intd{\intd t}(
(t|\x|)^{d-1}
\vtf \varphi_{\rho}(t\x)\cdot\x/|\x|
)\mid_{t=1}
\]
and
\[
\Div(\ssym{\vtf \varphi})(\x)
=
\Div((\ssym{\vtf \varphi})_{\rho})(\x)
=
\frac1{|\x|^{d-1}}
\frac\intd{\intd r}(
r^{d-1}
\varphi(r)
)(|\x|)
.
\]

Letting $\vt e=(1,0,\ldots,0)\in\mathbb{R}^d$ we have
\begin{align*}
\int_{\partial B(\vt 0,r)}
\Div(\vtf \varphi)
\intd\lmo{d-1}
&=
\sigma_{d-1}r^{d-1}
\int_{\SO{d}}
\Div(\vtf \varphi)(Rr\vt e)
\intd R\\
& =
\sigma_{d-1}r^{d-1}
\Div\left(
\int_{\SO{d}}
R^{-1}\vtf \varphi(R\cdot)
\intd R
\right)(r\vt e)
\\
&=
\sigma_{d-1}r^{d-1}
\Div(\ssym{\vtf \varphi})(r \vt e)
=
\sigma_{d-1}
\frac\intd{\intd r}(
r^{d-1}
\varphi(r)
)
,
\end{align*}
and we obtain
\begin{align*}
\int_{\R^d} f\Div(\vtf \varphi)
&=
\int_0^\infty
\int_{\partial B(\vt 0,r)}
f\Div(\vtf \varphi)
\intd\lmo{d-1}
\intd r
=
\sigma_{d-1}
\int_0^\infty
\rad f(r)
\frac\intd{\intd r}(
r^{d-1}\varphi(r)
)
\intd r
.
\end{align*}
Taking the supremum on both sides of the previous display over $\vtf \varphi\in C^1_\cpt(\{\vt x:a<|\vt x|<b\}, \R^d)$ with $\|\vtf \varphi\|_{L^\infty(\R^d)}\leq1$ concludes the proof,
using \(\| \varphi\|_{L^\infty((a,b))}\leq\|\vtf\varphi\|_{L^\infty(\mathbb{R}^d)}\leq1\), and that, conversely,
for any given $\psi\in C^1_\cpt((a,b))$ with $\|\psi\|_{L^\infty((0,\infty))}\leq1$, defining $\vtf \varphi(\vt x)=\psi(|\vt x|)\vt x/|\vt x|$ we have 
$\vtf \varphi\in C^1_\cpt(\{\vt x:a<|\vt x|<b\}, \R^d)$, $\|\vtf \varphi\|_{L^\infty(\R^d)}\leq1$ and $ \varphi=\psi$.
\end{proof}

\begin{corollary}
\label{varfradial}
\label{gradfdisintegration}
\label{gradfdisintegrationintegration}
Let $f\in L^1_\loc(\mathbb{R}^d)$ be radial and $0\leq a < b\leq\infty$.
Then $f\in\BV_\loc(\R^d\setminus\{\vt 0\})$ if and only if $\rad f\in\BV_\loc((0,\infty))$, in which case
\begin{equation}
\label{gradfdisintegrationformula}
\var_{\{\vt x:a<|\vt x|< b\}} (f)
=
\sigma_{d-1}
\int_{(a,b)}
r^{d-1}
\intd|\rad f'|(r)\,,
\end{equation}
where both terms may be infinite.
\end{corollary}

\begin{proof}
\zcref[S]{gradfdisintegrationformulaweak} directly implies the equivalence statement regarding $\BV_\loc$, and it order to prove \zcref{gradfdisintegrationformula} it suffices to show
\begin{equation}
\label{weak_definition_weighted}
\int_{(a,b)}
r^{d-1}
\intd|\rad f'|(r)
=
\sup\left\{
\int
\varphi'
\rad f
:
\varphi\in C^1_\cpt((a,b)),\ |\varphi(r)|\leq r^{d-1}
\right\}
.
\end{equation}
First consider $0<a < b<\infty$. Then $\varphi\in C^1_\cpt((a,b))$ if and only if $r\mapsto\varphi(r)/r^{d-1}\in C^1_\cpt((a,b))$, and
\[
\int(\varphi(r)/r^{d-1})r^{d-1}\intd{\rad f'}(r)
=
\int\varphi'\rad f
.
\]
Moreover, $r^{d-1}\rad f'(r)\chi_{(a,b)}(r)$ represents a Radon measure, call it $\nu$.
Since $C^1_\cpt((a,b))$ is dense in $L^1(\nu)$, there exist $\psi_1,\psi_2,\ldots\in C^1_\cpt((a,b))$ that approximate the Radon-Nikodym derivative $\alpha:=\frac{\intd \nu}{\intd |\nu|}$.
This implies
\(
\int
\psi_k
\intd \nu
\rightarrow
\|\nu\|
.
\)
Moreover, since $|\alpha(x)|=1$ for $\nu$-almost every $x$, the sequence $\{\psi_n\}_{n\geq 1}$ can be taken with \(\|\psi_n\|_{L^\infty((0,\infty))}\leq1\). Thus, taking suprema over $\varphi$ in the previous equality we obtain \zcref{weak_definition_weighted} which finishes the proof. Letting $a \to 0^+$ and $b \to \infty$ yields \zcref{weak_definition_weighted} in these cases as well (possibly with both sides being infinite).
\end{proof}

\subsection{Functions of bounded variation in one dimension}
\label{appendix_bvoned}

We follow \cite[Section~3.2]{AmbrosioFuscoPallara}.
For an open interval $I\subset\mathbb{R}$ recall the definition of variation \zcref{definition1dvariation}, which we denote by $\varalt_I(f)$ in this \zcref*[noref,nocap]{appendix_bvoned}.
Then every $f\in\BV(I)$ has a representative $\rep f$ such that $\var_J(f)=\varalt_J(\rep f)$ for all open intervals $J\subset I$.
Conversely, let $f:I\rightarrow\mathbb{R}$ with $\varalt_I(f)<\infty$.
Then $f\in\BV(I)$ with $\var_I(f)\leq\varalt_I(f)$, and for every $x\in I$ the one sided limits
\begin{align*}
f^-(x)
&=
\lim_{x\in I,\ y<x,\ y\rightarrow x}
f(y)
,&
f^+(x)
&=
\lim_{x\in I,\ y>x,\ y\rightarrow x}
f(y)
\end{align*}
exist.
Moreover, $\var_I(f)=\varalt_I(f)$ if and only if for all $x\in I$ exists a $0\leq \theta\leq 1$ with $f(x)=\theta f^-(x)+(1-\theta)f^+(x)$.

In particular, this means the function given by
\begin{align*}
\lowerrep f(x)
=
\min\{f^-(x),f^+(x)\}
\end{align*}
satisfies
\begin{align*}
\varalt_I(\lowerrep f)
\leq
\varalt_I(f)
.
\end{align*}
This fact can also be proven directly from the definition of $\varalt_I(f)$ without passing through $\var_I(f)$.
Note, that if $\varalt_I(f)<\infty$ then $f$ has at most countably many points of discontinuity, which means for a.e.\ $x\in I$ we have $\lowerrep f(x)=f(x)$, and thus $\M\lowerrep f=\M f$.
We record this statement in the following \zcref*[noref,nocap]{representativeassuminglowerlimit}:

\begin{lemma}
\label{representativeassuminglowerlimit}
For $f:\mathbb{R}\rightarrow\mathbb{R}$ with $\varalt(f)<\infty$ we have \(\varalt(\lowerrep f)\leq\varalt(f)\) and \(\M\lowerrep f=\M f\).
\end{lemma}

In the case of characteristic functions $f:\mathbb{R}\rightarrow\{0,1\}$ with $\var(f)<\infty$ we know from \zcref{Prop11_contbv} that \(\lowerrep{\M_\varphi f}=\M_\varphi f\) which in particular means both notions of variation agree for $\M_\varphi f$,
\begin{equation}
\label{Prop11_cont}
\varalt(\M_\varphi f)
=
\var(\M_\varphi f)
.
\end{equation}

Next, let $f\in\BV(\mathbb{R}^d)$ be radial.
Then by \zcref{gradfdisintegrationintegration} we have $\rad f\in\BV_\loc((0,\infty))$, and $\grad f$ does not depend on $\rad f(0)$.
That means we may consider $\rad f$ as a map from $(0,\infty)$ to $\mathbb{R}$.
Combining the previous arguments, for radial $g:\mathbb{R}^d\rightarrow\mathbb{R}$ given by $g(\vt 0)=f(\vt 0)$ and
\[
g(\vt x)
=
\rad g(|\vt x|)=\lowerrep{\big(\rep{(\rad f)}\big)}(|\vt x|)\,,
\]
we obtain $\rad g(x)=\rad f(x)$ for a.e.\ $x>0$ and hence $g(\vt x)=f(\vt x)$ for a.e.\ $\vt x\in\mathbb{R}^d$.
This means $g\in\BV(\mathbb{R}^d)$ with $\var(g)=\var(f)$ and $\M g=\M f$.
Moreover, $\rad g\in\BV_\loc([0,\infty))$ with $\rad g'=\rad f'$, for any open $I\subset[0,\infty)$ away from $0$ we have $\varalt_I(\rad g)=\var(\rad f)<\infty$, and
\(
\rad g
=
\lowerrep{(\rad g)}
.
\)
In the case of radial characteristic functions $f$ this means the following:

\begin{lemma}
\label{radialcharfrepresentative}
Let $f\in\BV(\mathbb{R}^d)$ be a radial characteristic function with $\var(f)<\infty$.
Then exist $N\in\N\cup\{\infty\}$ and sequences $\{a_n\}_{-N-1<n\leq0}$ and $\{b_n\}_{-N-1<n\leq0}$, with $0< a_n<b_n<a_{n+1}$ for each $n$, such that $g:\mathbb{R}^d\rightarrow\mathbb{R}$ given by
\[
\rad g = \sum_{n=-N}^0 \chi_{(a_n, b_n)}
\]
satisfies \(\var(g)=\var(f)\) and $\M f=\M g$.
\end{lemma}

\section{Continuity of the maximal function}

We now present a simple, yet useful, \zcref*[noref,nocap]{Prop11_contbv}. We state it for the signed maximal operator $\M^s_\varphi$ (recall definition \zcref{20261009_03:10}), but it plainly holds for the maximal operator $\M_\varphi$ as a particular case. 

\begin{proposition}
\label{Prop11_contbv}
Let $\varphi\in L^1(\mathbb{R})$ be non-negative and let $f\in\BV(\R)$.
Then, for every $z\in\mathbb{R}$, we have
\[
\M^s_\varphi f(z)
=
\liminf_{z\neq x\rightarrow z}
\M^s_\varphi f(x)
.
\]
\end{proposition}

\begin{proof}
It suffices to consider $z=0$.
In general, maximal functions of convolution type are lower semicontinuous, so we have
\[\M^s_{\varphi}f(0) \leq \liminf_{0\neq x \to 0} \M^s_{\varphi}f(x).\]
Then it suffices to show that 
\begin{equation}\label{20260508_15:08}
\liminf_{0\neq x \to 0} \M^s_{\varphi}f(x) \leq \M^s_{\varphi}f(0).
\end{equation}

Since $f\in\BV(\R)$ we have $C:=\|f\|_{L^\infty(\mathbb{R}^d)}<\infty$.
Thus, for any $r>0$ and $x,y\in\mathbb{R}$ we have
\begin{equation}
\label{bigr_converge}
\begin{split}
|
f*\varphi_r(x)
-
f*\varphi_r(y)
|
&=
\left|
\int
f(-z)[
\varphi_r(z+x)
-
\varphi_r(z+y)
]
\d z
\right|
\\
&\leq
C
\int
|
\varphi_r(z+x)
-
\varphi_r(z+y)
|
\d z
\\
&=
C
\int
|
\varphi(z)
-
\varphi(z+(y-x)/r)
|
\d y
\\
&\rightarrow0
\qquad
\text{as }
(y-x)/r
\rightarrow0
.
\end{split}
\end{equation}

Let $x_n \rightarrow0$ for which $\Ms_\varphi f(x_n)$ converges.
Then there exist $r_n>0$ with
\begin{equation}
\label{good_rn}
\lim_{n\rightarrow\infty} f*\varphi_{r_n}(x_n)
=
\lim_{n\rightarrow \infty} \M^s_{\varphi}f(x_n)
.
\end{equation}
We consider the following dichotomy:
\begin{enumerate}
\item
\label{largeradius}
There exist $R>0$ and $x_n\rightarrow0$ for which $\Ms_\varphi f(x_n)$ converges, and $r_n\geq R$ such that \zcref{good_rn} holds.
\item
\label{smallradius}
For all $x_n\rightarrow0$ for which $\Ms_\varphi f(x_n)$ converges, if $r_n>0$ satisfies \zcref{good_rn} then $r_n\rightarrow0$.
\end{enumerate}
Note that both statements also quantify over the choice $x_n=0$.

Assume \zcref{largeradius} holds.
Then for any other sequence $y_n\rightarrow0$ by \zcref{bigr_converge} we can conclude
\[
\limsup_{n\rightarrow\infty}
\M^s_{\varphi}f(y_n)
\geq
\limsup_{n\rightarrow\infty}
f*\varphi_{r_n}(y_n)
=
\limsup_{n\rightarrow\infty}
f*\varphi_{r_n}(x_n)
=
\lim_{n\rightarrow \infty} \Ms_{\varphi}f(x_n)
.
\]
Taking $y_n\neq0$ realizing $\liminf_{0\neq x \to 0} \M^s_{\varphi}f(x)$, i.e.
\[
\lim_{n\rightarrow\infty} \M^s_{\varphi}f(y_n)
=
\liminf_{0\neq x \to 0} \M^s_{\varphi}f(x)
,
\]
this implies that also $x_n$ realizes it.
Then taking $y_n=0$ this concludes the proof of \zcref{20260508_15:08} in this case.

It remains to consider \zcref{smallradius}.
Since $f\in\BV(\R)$ there exists $C\in\mathbb{R}$ such that $f=C\chi_{[0,\infty)}+g$ almost everywhere, where $g$ is continuous in $0$.
By symmetry it suffices to consider $C\geq0$.
Let $x_n\rightarrow0$ and $r_n\rightarrow0$.
Then
\begin{align*}
C\chi_{[0,\infty)}*\varphi_{r_n}(x_n)
&=
C\int_{-\infty}^{x_n/r_n}
\varphi
,&
g*\varphi_{r_n}(x_n)
&\rightarrow
g(0)
.
\end{align*}
In particular taking $x_n=0$ this means
\[
\Ms_\varphi f(0)
\geq
g(0)
+
C\int_{-\infty}^0\varphi
.
\]
There exists a sequence $0>x_n\rightarrow0$ for which $\Ms_\varphi f(x_n)$ converges and
$0<r_n\rightarrow0$ satisfying \zcref{good_rn}.
Since $\varphi\geq0$, we can conclude
\[
\liminf_{0\neq x\rightarrow0}
\Ms_\varphi f(x)
\leq
\lim_{n\rightarrow\infty}
\Ms_\varphi f(x_n)
=
\lim_{n\rightarrow\infty}
f*\varphi_{r_n}(x_n)
\leq
g(0)
+
C\int_{-\infty}^0\varphi
\leq
\Ms_\varphi f(0)
,
\]
which finishes the proof.
\end{proof}

\printbibliography

\end{document}